\magnification=\magstephalf
\input dfeo.mac
\input epsf
\let\noarrow = t
\input eplain
\indexproofingtrue
\newcount\refno\refno=1
\def\incr{\advance\refno by 1}
\edef\Berth{\number\refno}\incr
\edef\Bourb{\number\refno}\incr
\edef\deJO{\number\refno}\incr
\edef\Font{\number\refno}\incr
\edef\GorOo{\number\refno}\incr
\edef\Illus{\number\refno}\incr
\edef\KottIsoc{\number\refno}\incr
\edef\Kott{\number\refno}\incr
\edef\Kraft{\number\refno}\incr
\edef\GSAS{\number\refno}\incr
\edef\STPEL{\number\refno}\incr
\edef\FOAnn{\number\refno}\incr
\edef\FOTexel{\number\refno}\incr
\edef\FOTexelB{\number\refno}\incr
\edef\RapBourb{\number\refno}\incr
\edef\RR{\number\refno}\incr
\edef\Wedh{\number\refno}\incr

\centerline{{\sectitlefont A DIMENSION FORMULA FOR EKEDAHL-OORT STRATA}}
\Askip

\centerline{{\it by}\quad{\namefont Ben Moonen\footnote{*}{{\eightrm Research 
made possible by a fellowship of the Royal Netherlands Academy of Arts and 
Sciences}}\phantom{\quad{\it by}}}}
\Askip

To an abelian variety $X$ over a field~$k$ of characteristic~$p$ we can associate invariants such as the $p$-rank 
or the isogeny class of its $p$-divisible group. Such invariants can be used to define stratifications of the 
moduli space $\cA_{g}$ in characteristic~$p$. One example of such a stratification is the Newton stratification. 
Two points of $\cA_g(k)$ are in the same Newton stratum iff the associated $p$-divisible groups are isogenous 
over~$\kbar$. Mainly through the work of de Jong and Oort~[\deJO] and Oort~[\FOAnn], we now have a fairly complete 
picture of this stratification. For an overview we refer to Oort~[\FOTexelB] or Rapoport's Bourbaki 
lecture~[\RapBourb].

Of more recent date is the EO-stratification, after Ekedahl and Oort. The starting point is that, given $g = 
\dim(X)$ and working over $k = \kbar$ with $\charact(k) =p$, there is a finite list of isomorphism classes of 
group schemes~$X[p]$. This was proven by Kraft~[\Kraft] (unpublished), was seemingly forgotten for some time, and 
was then reobtained by Oort around~1995. The EO-stratification is defined by declaring that two points of 
$\cA_g(k)$ are in the same stratum iff the associated group schemes $X[p]$ are isomorphic over~$\kbar$. Several 
basic properties of this stratification were obtained by Ekedahl and Oort; see Oort~[\FOTexel]. Other questions, 
such as which strata occur in the boundary of a given one, remain open.

Naturally, people have tried to extend these theories to other Shimura varieties. A good class to start working on 
are the PEL moduli spaces at primes of good reduction. A generalized Newton stratification can be defined using 
the theory of isocrystals with additional structure, developed by Kottwitz~[\KottIsoc] and Rapoport and 
Richartz~[\RR]. The picture is, at present, still much less complete than in the Siegel modular case, though; see 
[\RapBourb],~\S5. 

Work on the generalized EO-stratification was taken up by Wedhorn~[\Wedh] and the author~[\GSAS]. (The Hilbert 
modular case was studied by Goren and Oort in~[\GorOo].) The data that 
are fixed in the moduli problem give rise to an algebraic group~$G$ and a conjugacy class~$\mX$ of parabolic 
subgroups. The ``$p$-kernel'' objects to consider are triples $\ul{Y} = (Y,\iota,\lambda)$ consisting of a group 
scheme~$Y$ killed by~$p$, equipped with an action~$\iota$ of a semi-simple $\Fp$-algebra and a polarization 
$\lambda\colon Y \isomarrow Y^D$. Writing $W_G$ for the Weyl group of~$G$ and $W_\mX \subset W_G$ for the subgroup 
associated to~$\mX$, the main result of~[\GSAS] is that, fixing suitable discrete invariants, such 
triples~$\ul{Y}$ are classified by the $W_\mX$-cosets in~$W_G$. Wedhorn showed in~[\Wedh] that this leads to a 
generalized EO-stratification on good reductions~$\cA_0$ of PEL moduli spaces. More precisely, $k$-valued points 
of~$\cA_0$ naturally give rise to triples~$\ul{Y}$ as above, and with the appropriate assumptions we arrive at a 
stratification $\cA_0 = \amalg_{w \in W_\mX\backslash W_G} \cA_0(w)$. Also in~[\Wedh] we find a result about the 
dimension of the strata: If $\ul{Y}_w$ is the triple corresponding to the coset $w \in W_\mX\backslash W_G$ and if 
$\cA_0(w) \neq \emptyset$, then every irreducible component of~$\cA_0(w)$ has codimension equal to 
$\dim\big(\bAut(\ul{Y}_w)\big)$ in~$\cA_0$.

The main purpose of the present paper is to give an explicit dimension formula for the strata. To state the 
result, let us recall that once we fix a generating set of reflections $S \subset W_G$, every coset $w \in 
W_\mX\backslash W_G$ has a distinguished representative~$\dot w$. The actual element $\dot w$ depends on the 
choice of~$S$, but its length $\ell(\dot w)$ does not. Our main result is then the following.
\medskip

\noindent
{\it Theorem. --- If $\cA_0(w) \neq \emptyset$ then its irreducible components all have dimension equal 
to~$\ell(\dot w)$.}
\medskip

The strength of this formula lies in the fact that the lengths $\ell(\dot w)$ are easily computable. In 
particular, an immediate consequence of our result is that there is at most one $0$-dimensional stratum and that 
there is a unique stratum of maximal dimension. The latter takes the role of a generalized ordinary stratum and 
plays a central role in our paper~[\STPEL], in which a generalization of Serre-Tate theory is developed.
\medskip

As suggested by the above, what we really compute in this paper are the dimensions of the automorphism group 
schemes $\bAut(\ul{Y})$. These computations are based on an explicit description of the Dieudonn\'e modules of 
triples~$\ul{Y}$ as in the above. The classification results of~[\GSAS] are recalled in \S1, the actual 
computation of the dimension of $\bAut(\ul{Y})$ is done in \S2. The application to the study of EO-stratifications 
is discussed in~\S3.
\medskip

\noindent
{\it Acknowledgements.\/} I thank F.~Oort, R.~Pink and T.~Wedhorn for stimulating discussions on the subject of 
this paper. The research for this paper was made possible by a Fellowship of the Royal Netherlands Academy of Arts 
and Sciences (KNAW). During my work on this paper I have been affiliated to the University of Utrecht (until June 
2001) and the University of Amsterdam (from July 2001). I thank these institutions for their support.

\section{$p$-Kernel group schemes with additional structures}{}

\ssection{Generalities on BT$_1$}{} 
\Bskip

\noindent
We fix a prime number~$p$. When dealing with polarized group schemes we assume~$p \neq 2$.

\sssection{BTdef}
Let $S$ be a scheme. Write $S_0 \subset S$ for the closed subscheme defined by the ideal $(p) \subset O_S$. By a 
{\it BT$_1$ over~$S$\/} (short for ``truncated Barsotti-Tate group of level~$1$'') we mean a commutative finite 
locally free $S$-group scheme~$Y$ such that, with the notation $Y_0 := Y \times_S S_0$, the sequence
$$
Y_0 \mapright{F} Y_0^{(p)} \mapright{V} Y_0
$$
is exact. Here $F$ and $V$ denote the relative Frobenius and Verschiebung of $Y_0$ over~$S_0$. For further details 
see Illusie~[\Illus].

If $Y$ is a BT$_1$ over~$S$ then we write $Y^D$ for its Cartier dual. There is a canonical isomorphism 
$Y \isomarrow Y^{DD}$, which we take as an identification. If $\epsilon \in \{\pm 1\}$ then by an {\it $\epsilon$-duality\/} of~$Y$ we mean an isomorphism $\lambda\colon Y \isomarrow Y^D$ such that $\lambda = \epsilon \cdot \lambda^D$. Such an $\epsilon$-duality induces an involution $f \mapsto f^\dagger$ on the ring $\End_S(Y)$. We also refer to an $\epsilon$-duality as a polarization.

\sssection{BTOstr} 
Let $B$ be an $\Fp$-algebra. By a {\it BT$_1$ with $B$-structure\/} over a basis~$S$ we mean a pair $\ul{Y} = 
(Y,\iota)$ where $Y$ is a BT$_1$ over~$S$ and $\iota\colon B \to \End_S(Y)$ is a homomorphism of $\Fp$-algebras.

Suppose $B$ is equipped with an $\Fp$-linear involution $b \mapsto b^\ast$. Let $\epsilon \in \{\pm 1\}$. By a {\it 
BT$_1$ with $(B,\ast,\epsilon)$-structure\/} over~$S$ we mean a triple $\ul{Y} = (Y,\iota,\lambda)$ where 
$(Y,\iota)$ is a BT$_1$ with $B$-structure and $\lambda\colon Y \to Y^D$ is an $\epsilon$-duality, such that 
$\iota(b^\ast) = \iota(b)^\dagger$ for all $b \in B$.

\sssection{DieuTh}
We use contravariant Dieudonn\'e theory as in Fontaine~[\Font]. Let $K$ be a perfect field, $\charact(K) = p$. Then 
a BT$_1$ with $B$-structure over~$K$ corresponds to a $4$-tuple $(N,F,V,\iota)$, where

\noindent
--- $N$ is a finite dimensional $K$-vector space,

\noindent
--- $F\colon N \to N$ is a $\Frob_K$-linear endomorphism, 

\noindent
--- $V \colon N \to N$ is a $\Frob_K^{-1}$-linear endomorphism, and

\noindent
--- $\iota\colon B \to \End(N,F,V)$ is an $\Fp$-linear homomorphism, 

\noindent
with $\Ker(F) = \Image(V)$ and $\Image(F) = \Ker(V)$. Using these last relations one can show that there exists a 
filtration
$$
\calC_\gdot\colon
\quad (0) = \calC_0 \subset \calC_1 \subset \calC_2 \subset \cdots \subset \calC_r = N
$$
that is the coarsest filtration with the properties that
\item{(\romno1)}for every $j$ there exists an index $f(j) \in \{0,1,\ldots,r\}$ with $F(\calC_j) = 
\calC_{f(j)}$;
\item{(\romno2)}for every $j$ there exists an index $v(j) \in \{0,1,\ldots,r\}$ with 
$V^{-1}(\calC_j) = \calC_{v(j)}$.

\noindent
We refer to this filtration as the {\it canonical filtration\/} of~$N$. See also [\GSAS],~2.5.

Similarly, a BT$_1$ with $(B,\ast,\epsilon)$-structure corresponds to a $5$-tuple $(N,F,V,\phi,\iota)$, where 
$(N,F,V,\iota)$ is as above and where $\phi\colon N \times N \rightarrow K$ is a perfect, $\epsilon$-symmetric 
bilinear form, such that 
$$
\leqalignno{
\phi(Fn_1,n_2) = \phi(n_1,Vn_2)^p\qquad
&\hbox{{\rm for all~$n_1$,~$n_2 \in N$;}}&(\refn{DieuTh}.1)\cr
\phi(bn_1,n_2) = \phi(n_1,b^\ast n_2)\qquad
&\hbox{{\rm for all $b\in B$ and $n_1$, $n_2 \in N$.}}&(\refn{DieuTh}.2)\cr}
$$

\ssection{BT$_1$ with given endomorphisms}{BT1Endo} 
\Bskip

\noindent
Let $B$ be a finite dimensional semi-simple $\Fp$-algebra. Let $k$ be an algebraically closed field of characteristic~$p$. The first problem studied in~[\GSAS] is the classification of BT$_1$ with $B$-structure over~$k$. This generalizes the work of Kraft~[\Kraft], who classified group schemes killed by~$p$ without additional structure. We shall briefly review our results. 

\sssection{pairs}
Write $\kappa$ for the center of~$B$. Then $\kappa$ is a product of finite fields, say $\kappa = \kappa_1 \times \cdots \times \kappa_\nu$. Let~$\cI = \cI_1 \cup \cdots \cup\cI_\nu$ be the set of homomorphisms $\kappa \rightarrow k$.

Consider pairs $(N,L)$ consisting of a finitely generated $B \otimes_\Fp k$-module~$N$ and a submodule $L \subset N$. Note that the simple factors of $B \otimes_\Fp k$ are indexed by~$\cI$, so we get canonical decompositions $N = \oplus_{i\in\cI} N_i$ and $L = \oplus_{i\in\cI} L_i$. Define two functions $d$, $\gf \colon \cI \to \mZ_{\geq 0}$ by $d(i) = \len(N_i)$ and $\gf(i) = \len(L_i)$, takings lengths as $B \otimes_{\Fp} k$-modules. The pair $(d,\gf)$ determines the pair $(N,L)$ up to isomorphism.

To the pair $(N,L)$ we associate an algebraic group~$G$ over~$k$ and a conjugacy class~$\mX$ of parabolic subgroups of~$G$. First we define
$$
G := \GL_{B \otimes_\Fp k}(N) \cong \dirprod_{i\in\cI} \GL_{d(i),k}\, .
$$
Then the stabilizer $P := \Stab(L)$ is a parabolic subgroup of~$G$, and we define $\mX$ as the conjugacy class of 
parabolic subgroups of~$G$ containing $P$.

\sssection{BT1pair}
Let $\ul{Y} = (Y,\iota)$ be a BT$_1$ with $B$-structure over~$k$. Write~$N$ for the Dieudonn\'e module of~$Y$ and 
let $L := \Ker(F) \subset N$. Let $(d,\gf)$ be the corresponding pair of functions. It can be shown (see [\GSAS], 4.3) that the function~$d$ is constant on each of the subsets $\cI_n \subset \cI$. We refer to $(d,\gf)$ as the {\it type\/} of~$(Y,\iota)$.

\sssection{wdefGE}
Fix a pair $(d,\gf)$ with $d$ constant on each subset $\cI_n \subset \cI$. Fix a pair of $B \otimes_\Fp k$-modules 
$L_0 \subset N_0$ of type~$(d,\gf)$. Let $(G,\mX)$ be the associated algebraic group and conjugacy class of 
parabolic subgroups. Let $W_G$ be the Weyl group of~$G$, and let $W_\mX \subset W_G$ be the subgroup corresponding 
to~$\mX$.

To a pair $\ul{Y} = (Y,\iota)$ of type $(d,\gf)$ we associate an element $w(\ul{Y}) \in W_\mX \backslash W_G$. This 
is done as follows. Write~$N$ for the Dieudonn\'e module of~$Y$ and let $L := \Ker(F) \subset N$. Choose an 
isomorphism $\xi\colon N \isomarrow N_0$ that restricts to~$L \isomarrow L_0$. This allows us to view the canonical 
filtration~$\calC_\gdot$ of~$N$ as a filtration of~$N_0$. Choose any refinement~$\cF_\gdot$ of~$\calC_\gdot$ to a 
complete flag. The relative position of~$L_0$ and $\cF_\gdot$ is given by an element $w(L_0,\cF_\gdot) \in W_\mX 
\backslash W_G$. It can be shown that this element is independent of the choice of~$\xi$ and the 
refinement~$\cF_\gdot$; see [\GSAS], especially~4.6 for details. Now define $w(\ul{Y}) := w(L_0,\cF_\gdot)$.

With these notations, the first main result of [\GSAS] can be stated as follows.

\sssection{GEThm}
{\it Theorem. --- Assume that $k = \kbar$. Fix a type $(d,\gf)$. The map $\ul{Y} \mapsto w(\ul{Y})$ gives a bijection
$$
\left\{
\vcenter{
\setbox0=\hbox{{\rm isomorphism classes of}}
\copy0
\hbox to \wd0{{\rm \hfil $\ul{Y}$ of type $(d,\gf)$\hfil}}}
\right\}
\longisomarrow
W_\mX \backslash W_G\, .
$$}
\vskip-\belowdisplayskip

\ssection{BT$_1$ with given endomorphisms and a polarization}{BT1EP} 
\Bskip

\noindent 
Let $B$ be a finite dimensional semi-simple $\Fp$-algebra equipped with an involution $b \mapsto b^\ast$. Let 
$\epsilon \in \{\pm 1\}$. Let $k$ be an algebraically closed field of characteristic~$p > 2$. The second problem 
studied in~[\GSAS] is the classification of BT$_1$ with $(B,\ast,\epsilon)$-structure over~$k$. We shall briefly 
review the result. As in section~\refn{BT1Endo}, this involves an algebraic group~$G$. In the polarized case there 
are two versions of the result: one using a possibly non-connected group~$G$, the other using the identity 
component~$G^0$. We shall need both variants; see \refn{delta}--\refn{ww0} below; see also~[\GSAS], especially~3.8 
and~5.7.

\sssection{CDAfact}
Let $(B,\ast,\epsilon)$ be as above. We can decompose $(B,\ast)$ as a product of simple factors, say
$$
(B,\ast) = (B_1,\ast_1) \times \cdots \times (B_\nu,\ast_\nu)\, .
$$
If $\ast_n$ is an orthogonal involution, set $\epsilon_n := +1$; if $\ast_n$ is symplectic, set $\epsilon_n := -1$. 
The simple factors $(B_n,\ast_n)$ come in three kinds:
\smallskip

\typesitem{Type C:} $B_n \cong M_r(\kappa_n)$ for some finite field~$\kappa_n$, with $\ast_n$ of the first kind and 
$\epsilon\cdot \epsilon_n = -1$;

\typesitem{Type D:} $B_n \cong M_r(\kappa_n)$ for some finite field~$\kappa_n$, with $\ast_n$ of the first kind and 
$\epsilon\cdot \epsilon_n = +1$;

\typesitem{Type A:} $B_n \cong M_r(\tilde\kappa_n)$, where $\tilde\kappa_n$ is an \'etale quadratic extension of a 
finite field~$\kappa_n$, and $\ast_n$ is of the second kind.

\noindent
Note that the labelling depends on~$\epsilon$. For factors of type~A, if $\kappa_n \cong \mF_q$ then either 
$\tilde\kappa_n \cong \mF_{q^2}$ or $\tilde\kappa_n = \kappa_n\times\kappa_n$.

Let $\tilde \kappa$ be the centre of~$B$, and define $\kappa := \{z\in \tilde\kappa \mid z^\ast = z\}$. We have 
$\tilde\kappa = \tilde\kappa_1 \times \cdots \times \tilde\kappa_\nu$ and $\kappa = \kappa_1 \times \cdots \times 
\kappa_\nu$, where the $\kappa_n$ are finite fields, $\tilde\kappa_n = \kappa_n$ if $(B_n,\ast_n)$ is of type~C 
or~D, and $\tilde\kappa_n$ is an \'etale quadratic extension of~$\kappa_n$ if $(B_n,\ast_n)$ is of type~A. Let $\cI 
= \cI_1 \cup \cdots \cup \cI_\nu$ be the set of homomorphisms $\kappa \rightarrow k$. For $X \in \{{\rm A}, {\rm 
C}, {\rm D}\}$ we say that $i \in \cI$ is of type~$X$ if $i \in \cI_n \subset \cI$ and $(B_n,\ast_n)$ is of 
type~$X$. Let $\cI^X \subset \cI$ be the subset of elements of type~$X$. Let~$\cItil$ be the set of homomorphisms 
$\tilde\kappa \rightarrow k$. We have a restriction map $\res\colon \cItil \rightarrow \cI$. For $\tau \in\cItil$ 
define $\bar\tau := \tau \circ \ast$. If $i\in\cI$ is of type~C or~D then there is a unique $\tau \in \cItil$ with 
$\res(\tau) = i$, and $\tau = \bar\tau$; if $i$~is of type~A then there are precisely two elements $\tau$, 
$\bar\tau \in \cItil$ that restrict to the embedding $i$ on~$\kappa$.

\sssection{pairs2}
Consider triples $(N,L,\phi)$ consisting of a finitely generated $B \otimes_\Fp k$-module~$N$, a perfect, 
$\epsilon$-symmetric bilinear form $\phi\colon N \times N \to k$ satisfying (\refn{DieuTh}.2), and a maximal 
isotropic submodule $L \subset N$. With a similar construction as in~\refn{pairs}, such a triple is classified, up to isomorphism, by a pair $(d,\gf)$ consisting of functions $d\colon \cI \to \mZ_{\geq 
0}$ and $\gf\colon \cItil \to \mZ_{\geq 0}$ such that $\gf(\tau) + \gf(\bar\tau) = d(i)$ for all $\tau\in\cItil$ and $i = \res(\tau) \in \cI$. Note that if $i$ is of type C or~D then there is a unique $\tau = \bar\tau$ with $\res(\tau) = i$ and we get the relation $d(i) = 2 \gf(\tau)$.

To a triple $(N,L,\phi)$ as above we associate a pair $(G,\mX)$. First, define $G := \UU_{B \otimes_\Fp 
k}(N,\phi)$, the algebraic group (over~$k$) of $B \otimes_\Fp k$-linear automorphisms of~$N$ that preserve the 
form~$\phi$. We have $G = \prod_{i\in\cI} G_i$, with $G_i$ isomorphic to $\Sp_{d(i),k}$ if $i$ is of type~C, to $\OO_{d(i),k}$ if $i$ is of type~D and to $\GL_{d(i),k}$ if $i$ is of type~A. Note that $d(i)$ is even if $i$~is of type~C or~D. In the presence of (non-connected) orthogonal factors, the set~$\mX$ that we want to consider is not simply a conjugacy class of parabolic subgroups of~$G$. Instead we consider (partial) hermitian flags in~$N$, i.e., filtrations by $B \otimes_\Fp k$-submodules
$$
(0) = C_0 \subset C_1 \subset \cdots \subset C_r = N
$$
with the property that $C_j^\perp = C_{r-j}$ for all~$j$. (Here $C_j^\perp := \{n\in N\mid \phi(C_j,n) = 0\}$.) The 
set $\Flag(N,\phi)$ of all such hermitian flags has the structure of a (generally non-thick) building, on which $G$ 
acts in a strongly transitive, type-preserving manner. Then we define $\mX$ to be the $G$-orbit in $\Flag(N,\phi)$ 
of the flag $(0) \subset L \subset N$. Equivalently, $\mX$ is the set of all maximal isotropic subspaces in~$N$. 
Note that if there are no factors of type~D then $C_\gdot \mapsto \Stab(C_\gdot)$ gives a bijective correspondence 
between hermitian flags and parabolic subgroups of~$G$; in this case we could also define~$\mX$ to be the conjugacy 
class of parabolic subgroups of~$G$ that contains $\Stab(L)$.

\sssection{BT1pair2}
Let $\ul{Y}$ be a BT$_1$ with $(B,\ast,\epsilon)$-structure over~$k$. Let $N$ be the Dieudonn\'e module of~$Y$, let 
$L = N[F] := \Ker(F)$, and let $\phi$ be the $\epsilon$-hermitian form on~$N$ corresponding to the given 
$\epsilon$-duality. Then $(N,L,\phi)$ is a triple as in~\refn{pairs2}. Let $(d,\gf)$ be the corresponding pair of 
functions; we refer to this pair as the type of~$\ul{Y}$. It can be shown ([\GSAS], 4.3, 5.3 and 6.5) that the 
function~$d$ is constant on each of the subsets $\cI_n \subset \cI$.

\sssection{wdefGPE}
Fix $(d,\gf)$ as in~\refn{pairs2}, with $d$ constant on each subset $\cI_n \subset \cI$. Choose a 
corresponding triple $(N_0,L_0,\phi_0)$. To this triple we associate a pair $(G,\mX)$ as explained above. The 
Coxeter group associated to the building $\Flag(N_0,\phi_0)$ is just the Weyl group $W_G$ of~$G$. Note, however, 
that we work with a possibly non-connected group~$G$; an orthogonal factor $\OO_{2q}$ contributes a factor of type 
B$_q$ (not D$_q$) to~$W_G$. Write $W_\mX \subset W_G$ for the subgroup corresponding to~$\mX$.

Let $\ul{Y}$ be a BT$_1$ with $(B,\ast,\epsilon)$-structure over~$k$ of type $(d,\gf)$. Let $(N,L,\phi)$ be the 
associated triple, as in~\refn{BT1pair2}. To $\ul{Y}$ we associate an element $w(\ul{Y}) \in W_\mX \backslash W_G$. 
This works essentially the same as in the non-polarized case. First we choose an isometry $\xi \colon (N,\phi) 
\isomarrow (N_0,\phi_0)$ that restricts to $L \isomarrow L_0$. Via $\xi$ we can view the canonical filtration 
$\calC_\gdot$ as an element of $\Flag(N_0,\phi_0)$. Choose a refinement of $\calC_\gdot$ to a full hermitian flag 
$\cF_\gdot$ in~$N_0$. (I.e., choose a chamber of which $\calC_\gdot$ is a face.) The relative position of~$L_0$ and 
$\cF_\gdot$ is measured by an element $w(L_0,\cF_\gdot) \in W_\mX \backslash W_G$. It can be shown that this 
element is independent of the choice of~$\xi$ and the refinement~$\cF_\gdot$; see [\GSAS], especially~5.6 and~6.6 
for details. 
Now define $w(\ul{Y}) := w(L_0,\cF_\gdot)$.

With these notations, the second main result of [\GSAS] is the following.

\sssection{GPEThm}
{\it Theorem. --- Let $k$ be an algebraically closed field, $\charact(k) > 2$. Fix a type $(d,\gf)$. Sending a 
BT$_1$ with $(B,\ast,\epsilon)$-structure $\ul{Y}$ to the element $w(\ul{Y})$ gives a bijection
$$
\left\{
\vcenter{
\setbox0=\hbox{{\rm isomorphism classes of}}
\copy0
\hbox to \wd0{{\rm \hfil $\ul{Y}$ of type $(d,\gf)$\hfil}}}
\right\}
\longisomarrow
W_\mX \backslash W_G\, .
$$}

(In [\GSAS] this result is only stated for BT$_1$ with $(B,\ast,-1)$-structure. However, by Morita 
equivalence one can always reduce to the case that $\epsilon = -1$.)

\sssection{delta}
We describe a variant of the theorem that uses the Weyl group of the identity component~$G^0$. To explain this we 
need to introduce another invariant, which is a function $\delta\colon \cI^{{\rm D}} \to \mZ/2\mZ$. If there are no factors of type~D then $\cI^{{\rm D}} = \emptyset$ and the invariant~$\delta$ is void. If $i \in \cI^{{\rm D}}$, let $\tau \in \cItil$ be the unique element with $\res(\tau) = i$. As we have seen, the Dieudonn\'e module~$N$ decomposes as $N = \oplus_{\tau \in \cItil} N_\tau$. Let $N_\tau[F] := \Ker(F_{|N_\tau})$ and $N_\tau[V] = \Ker(V_{|N_\tau})$. Now define $\delta(i)$ to be the length of the $B \otimes_\Fp k$-module $N_\tau[F]/\big(N_\tau[F] \cap N_\tau[V]\big)$ modulo~$2$.

As in the above, let $G := \UU_{B \otimes_\Fp k}(N_0,\phi_0)$. Define $\mX^0$ to be the conjugacy class of 
parabolic subgroups of~$G^0$ containing $\Stab(L_0)$. Write $W_{G^0}$ for the Weyl group of~$G^0$, and let $W_{\mX^0} \subset W_{G^0}$ be the subgroup corresponding to~$\mX^0$.

In addition to the type $(d,\gf)$, also fix $\delta\colon \cI^{{\rm D}} \to \mZ/2\mZ$. Let $\ul{Y}$ be a BT$_1$ 
with $(B,\ast,\epsilon)$-structure of type $(d,\gf,\delta)$. With a similar procedure as above we associate to~$\ul{Y}$ an element $w^0(\ul{Y}) \in W_{\mX^0} \backslash W_{G^0}$. (If there are no factors of type~D then $w^0(\ul{Y})$ is just the same as~$w(\ul{Y})$.) 

With these notations we have the following variant of Thm.~\refn{GPEThm}.

\sssection{GPEThm2}
{\it Theorem. --- Let $k = \kbar$, $\charact(k) > 2$. Fix a type $(d,\gf,\delta)$. Then the map $\ul{Y} \mapsto w^0(\ul{Y})$ gives a bijection
$$
\left\{
\vcenter{
\setbox0=\hbox{{\rm isomorphism classes of}}
\copy0
\hbox to \wd0{{\rm \hfil $\ul{Y}$ of type $(d,\gf,\delta)$\hfil}}}
\right\}
\longisomarrow
W_{\mX^0} \backslash W_{G^0}\, .
$$}

\sssection{ww0}
If $m = \# \cI^{{\rm D}}$ then $\# (W_\mX\backslash W_G) = 2^m \cdot \#(W_{\mX^0}\backslash W_{G^0})$. This 
corresponds to the fact that, given $(d,\gf)$, there are $2^m$ choices for~$\delta$.

In some applications it is more convenient to work with the element $w(\ul{Y})$; in other cases $w^0(\ul{Y})$ is 
easier to use. At any rate, once we fix $\delta$ we have a bijection
$$
\left\{ 
\vcenter{
\setbox0=\hbox{elements of type $\delta$}
\copy0
\hbox to\wd0{\hfil in $W_\mX\backslash W_G$\hfil}
}\right\}
\isomarrow W_{\mX^0}\backslash W_{G^0}
$$
that can be made completely explicit. For details we refer to [\GSAS],~3.8.
\section{Automorphism group schemes}{AutGrSch}

\ssection{Statement of the main result}{AutState}

\sssection{DistRepr}
Let $(W,S)$ be a Coxeter system, $\mX \subset S$ a subset, $W_\mX \subset W$ the subgroup generated by~$\mX$. Write $\ell$ for the length function on~$W$. Every coset $w \in W_\mX\backslash W$ has a unique representative $\dot w$ of minimal length, called the $(\mX,\emptyset)$-reduced representative. (See Bourbaki~[\Bourb], Chap.~\Romno 4, Exercise~3.) If $W_\mX \subset W$ and $w \in W_\mX\backslash W$ are given then $\dot w$ in general depends on the choice of~$S$, but $\ell(w) := \ell(\dot w)$ is independent of~$S$.
\Cskip

This section is devoted to the proof of the following result.

\sssection{AutDimThm}
{\it Theorem. --- {\rm (\romno1)}\enspace Let $B$ be a finite dimensional semi-simple $\Fp$-algebra. Let~$\ul{Y}$ be a BT$_1$ with $B$-structure over a field $k = \kbar$ with $\charact(k) = p$. Let $(d,\gf)$ be the type of $\ul{Y}$, define~$G$ and~$\mX$ as in\/ {\rm \refn{pairs}}, and let $w := w(\ul{Y}) \in W_\mX\backslash W_G$. Then the automorphism group scheme $\bAut(\ul{Y})$ has dimension equal to $\dim(\mX) - \ell(\dot w)$.

{\rm (\romno2)}\enspace Assume that $p > 2$. Let $\ast$ be an involution on~$B$ and let $\epsilon \in \{\pm 1\}$. Let~$\ul{Y}$ be a BT$_1$ with $(B,\ast,\epsilon)$-structure over~$k$. Let $(d,\gf,\delta)$ be the type of $\ul{Y}$, define~$G^0$ and~$\mX^0$ as in\/ {\rm \refn{delta}}, and let $w^0 := w^0(\ul{Y}) \in W_{\mX^0}\backslash W_{G^0}$. Then the automorphism group scheme $\bAut(\ul{Y})$ has dimension equal to $\dim(\mX^0) - \ell(\dot w^0)$.}

\sssection{DSpace}
Let $S$ be a scheme of characteristic~$p$. Following Wedhorn~[\Wedh] we define a {\it Dieu\-don\-n\'e space over~$S$\/} to be a triple $(N,F^\sharp,V^\flat)$ where $N$ is a locally free $\cO_S$-module of finite type and where $F^\sharp\colon N^{(p)} \rightarrow N$ and $V^\flat\colon N \rightarrow N^{(p)}$ are $\cO_S$-linear homomorphisms such that $F^\sharp \circ V^\flat = 0$ and $V^\flat \circ F^\sharp = 0$. (We use $F$ and $V$ for semi-linear endomorphisms, $F^\sharp$ and $V^\flat$ for their linearizations.)

In the following, $\cD$ will either denote a finite dimensional semi-simple $\Fp$-algebra~$B$ or a triple $(B,\ast,\epsilon)$ as in section~\refn{BT1EP}. We have the notion of a Dieudonn\'e space with $\cD$-structure over~$S$; see [\Wedh], \S~5 for details.

Let $\ul{Y}$ be a BT$_1$ with $\cD$-structure over~$k$. The Dieudonn\'e module of $\ul{Y}$ defines, by linearization of $F$ and~$V$, a Dieudonn\'e space with $\cD$-structure $\ul{N}$ over~$k$. Both $\ul{Y}$ and~$\ul{N}$ give rise to an automorphism group scheme. Concretely, if $f\colon T \rightarrow \Spec(k)$ is a $k$-scheme then the $T$-valued points of $\bAut(\ul{Y})$ are the automorphisms of the pull-back $f^\ast(\ul{Y})$ as a BT$_1$ with $\cD$-structure over~$T$. Similarly, the $T$-valued points of $\bAut(\ul{N})$ are the automorphisms of the pull-back $f^\ast(\ul{N})$ as a Dieudonn\'e space with $\cD$-structure over~$T$.
\Cskip

We shall make use of the following result of Wedhorn.

\sssection{dim=dim}
{\it Proposition\/ {\rm (Wedhorn, [\Wedh])}. --- Notations as above. There is a natural homomorphism of $k$-group schemes
$$
\Delta \colon \bAut(\ul{Y}) \longrightarrow \bAut(\ul N)\, .
$$
If $K$ is any perfect field then $\Delta$ induces an isomorphism on $K$-valued points. This implies that $\dim\big(\bAut(\ul{Y})\big) = \dim\big(\bAut(\ul N)\big)$.}
\Cskip

Strictly speaking, [\Wedh] only deals with polarized BT$_1$. In the non-polarized case the same arguments apply. Alternatively, one can reduce to the polarized case by passing from a BT$_1$ with $B$-structure~$\ul{Y}$ to $Y \times Y^D$ with its natural $B \times B$-structure and the obvious $\epsilon$-duality.

To prove Thm.~\refn{AutDimThm} we shall compute the dimension of $\bAut(\ul N)$. Let us note that $\bAut(\ul{Y})$ and $\bAut(\ul{N})$ are in general non-reduced (see also \refn{NonRedRem}), and that $\Delta$ is in general not an isomorphism.  

\ssection{The non-polarized case}{Autnonpol}

\sssection{B=kappa}
We first prove (\romno1) of Thm.~\refn{AutDimThm}. The notions involved in the proof are illustrated in an example in~\refn{AutDimExa} below. There is an easy reduction to the case that $B = \kappa$ is a finite field. In fact, as the Brauer group of a finite field is trivial we have $B \cong M_{r_1}(\kappa_1) \times \cdots \times M_{r_\nu}(\kappa_\nu)$ where the $\kappa_n$ are finite fields. Fixing such an isomorphism, every BT$_1$ with $B$-structure decomposes as a product $\ul{Y} = (\ul{Y}_1)^{r_1} \times \cdots \times (\ul{Y}_\nu)^{r_\nu}$, where $\ul{Y}_n$ is a BT$_1$ with $\kappa_n$-structure. If the theorem holds for each~$\ul{Y}_n$ then one readily checks that it also holds for~$\ul{Y}$.

\sssection{l(w)npol}
Suppose $B = \kappa$ is a field of $p^m$ elements. Recall that we write $\cI$ for the set of embeddings $\kappa \to k$. Then $\cI$ is a set of $m$ elements that comes equipped with a natural cyclic ordering: if $i\in\cI$ then we write $i +1$ for $\Frob_k \circ i$. 

A type $(d,\gf)$ consists of a positive integer~$d$ and a function $\gf\colon \cI \to \{0,\ldots,d\}$. As in~\refn{pairs} we associate to the type $(d,\gf)$ a pair~$(G,\mX)$; in this case we have $G = \prod_{i \in \cI} \GL_{d,k}$. After a suitable choice of coordinates we can 
identify
$$
W_G = \dirprod_{i\in\cI} \gS_d
\qquad\hbox{and}\qquad
W_\mX = \dirprod_{i\in\cI} \gS\big\{1,\ldots,\gf(i)\big\} \times \gS\big\{\gf(i)+1,\ldots,d\big\}\, .
$$
In each factor $\gS_d$ we take the transpositions $(j\; j+1)$ as a set of generators.

As recalled in~\refn{DistRepr}, each coset $w \in W_\mX \backslash W_G$ has a distinguished representative $\dot w = (\dot w_i)_{i\in\cI}$ in~$W_G$. The permutations $\dot w_i$ that arise in this way are characterized by their property that
$$
j^\prime < j\quad \hbox{and}\quad \dot w_i(j^\prime) > \dot w_i(j) \qquad\Rightarrow\qquad \dot w_i(j) \leq \gf(i) < \dot w_i(j^\prime)\, .\leqno(\refn{l(w)npol}.1)
$$
The length of the element $\dot w$ is given by $\ell(\dot w) = \sum_{i\in\cI} \ell(\dot w_i)$, with
$$
\ell(\dot w_i) = \sum_{n=1}^{\gf(i)} \big(\dot w_i^{-1}(n) - n\big) = \sum_{n=\gf(i)+1}^d \big(n - \dot w_i^{-1}(n)\big)\, .\leqno(\refn{l(w)npol}.2)
$$

\sssection{SOnpol}
Let $w = w(\ul{Y}) \in W_\mX\backslash W_G$. Let $\ul{N} = (N,F,V,\iota)$ be the Dieudonn\'e module of~$\ul{Y}$. We recall from [\GSAS],~4.9 an explicit description of~$\ul{N}$ in terms of the distinguished representative~$\dot w$. Namely, let $N$ be the $k$-vector space with basis $e_{i,j}$ for $i\in\cI$ and $j\in\{1,\ldots,d\}$. Write $N_i :=\sum_{j=1}^d k\cdot e_{i,j}$, so that $N = \oplus_{i\in\cI} N_i$. We let $a \in \kappa$ act on $N_i$ as multiplication by~$i(a)$. Next define $\Frob_k$-linear maps $F_i\colon N_i \rightarrow N_{i+1}$ and $\Frob_k^{-1}$-linear maps $V_i\colon N_i \leftarrow N_{i+1}$ by
$$
F_i(e_{i,j}) = \cases{
0&if $\dot w_i(j) \leq \gf(i)$;\cr
e_{i+1,n}&if $\dot w_i(j) = \gf(i)+n$;}
\quad\hbox{and}\quad
V_i(e_{i+1,j}) = \cases{
0&if $j \leq d-\gf(i)$;\cr
e_{i,n}&if $j = d-\gf(i) + \dot w_i(n)$.}
$$
This defines a Dieudonn\'e module with $\kappa$-action such that the corresponding BT$_1$ is isomorphic to~$\ul{Y}$.

\sssection{picture}
We find it useful to draw pictures of the Dieudonn\'e modules $\ul{N}$ obtained in this way; see Figure~\refn{Fig:Exa}. The boxes represent base vectors as indicated. An edge connecting boxes $(i,j)$ and $(i+1,j^\prime)$ either represents the relation $F(e_{i,j}) = e_{i+1,j^\prime}$ (an $F$-edge) or the relation $V(e_{i+1,j^\prime}) = e_{i,j}$ (a $V$-edge). To determine in which of the two cases we are we represent the base vectors $e_{i,j}$ with $F(e_{i,j}) = 0$ by shaded boxes; thus in the module $N_i$ there are $\gf(i)$ shaded boxes. Note that in most cases we do not need the shading: (\refn{l(w)npol}.1) implies that two edges of the same kind never cross. Since the image of $F\colon N_i \to N_{i+1}$ is 
spanned by the {\it first\/} $d-\gf(i)$ base vectors, it follows that an $F$-edge in the picture always has slope $\leq 0$ and a $V$-edge always has slope $\geq 0$.
\midinsert
\line{\hfil\epsfbox{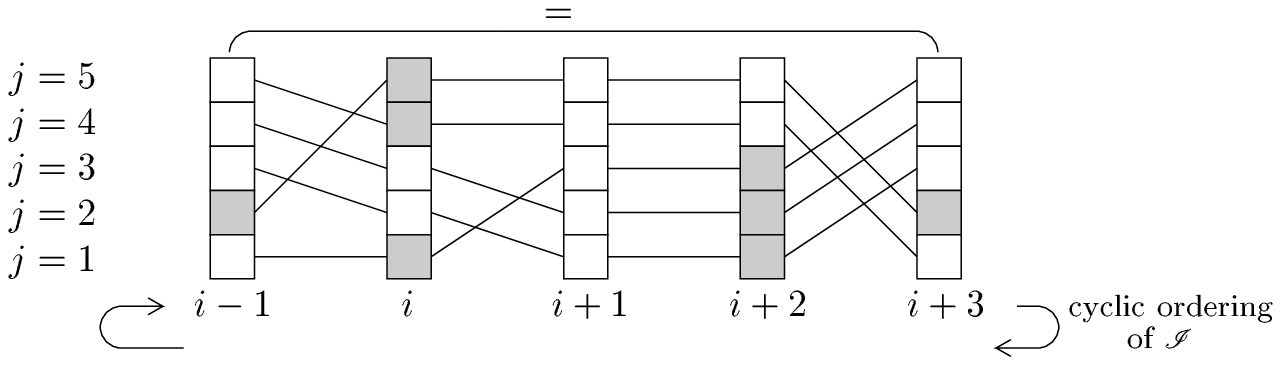}}
\figlabel{Fig:Exa}
\endinsert
 
In the example drawn in Figure~\refn{Fig:Exa} we take $d=5$ and $m=4$, with $\gf(i-1)=1$, $\gf(i)=3$, $\gf(i+1)=0$ and 
$\gf(i+2) = 3$. The distinguished representative~$\dot w$ is given by
$$
\dot w_{i-1} = \left[\matrix{1&2&3&4&5\cr 2&1&3&4&5\cr}\right]\, ,\quad
\dot w_{i} = \left[\matrix{1&2&3&4&5\cr 1&4&5&2&3\cr}\right]\, ,\quad
\dot w_{i+1} = \dot w_{i+2} = \left[\matrix{1&2&3&4&5\cr 1&2&3&4&5\cr}\right]\, .
$$
One should think of $\ul{N}$ as a roundabout (carousel); note that in Figure~\refn{Fig:Exa} (with $\#\cI = 4$) the summands $N_{i-1}$ and~$N_{i+3}$ have to be identified.

\sssection{canblocks} 
We recall some facts proven in [\GSAS], Lemma~4.5. The canonical filtration of~$N = \oplus_{i\in\cI} N_i$ induces filtrations 
$$
C_{i,\gdot}:\quad (0) = C_{i,0} \subset C_{i,1} \subset \cdots \subset C_{i,\ell} = N_i\, .
$$
In the description of $N$ given in \refn{SOnpol}, if $n = \dim(C_{i,r})$ then $C_{i,r} = k \cdot e_{i,1} + \cdots + k \cdot e_{i,n}$. The length~$\ell$ of the filtrations $C_{i,\gdot}$ is independent of~$i \in \cI$. We refer to the vector spaces $B_{i,j} := C_{i,j}/C_{i,j-1}$ as the {\it canonical blocks\/} of~$N$; they are indexed by the set $\cA := \cI \times \{1,\ldots,\ell\}$. For each pair $(i,j) \in \cA$ there is a unique index $\rho_i(j) \in \{1,\ldots,\ell\}$ such that
\line{\indent either: $F_i \colon N_i \to N_{i+1}$ induces a $\Frob_k$-linear bijection $B_{i,j} \isomarrow B_{i+1,\rho_i(j)}$,\hfil}
\line{\setbox0=\hbox{either:}\indent\hbox to\wd0{\hfil or:} $V_i \colon N_i \leftarrow N_{i+1}$ induces a $\Frob_k^{-1}$-linear bijection $B_{i,j} \leftisomarrow B_{i+1,\rho_i(j)}$.\hfil}

\noindent
Then $\rho(i,j) = \big(i+1,\rho_i(j)\big)$ defines a permutation~$\rho$ of the set~$\cA$, and for each $a = (i,j) \in \cA$ we have a $\Frob_k$-linear bijection $t_a\colon B_a \isomarrow B_{\rho(a)}$ which is either induced by~$F_i$ or is the inverse of the bijection induced by~$V_i$.

If $f \in \bAut(\ul{N})(R)$ for some $k$-algebra $R$ then $f$ preserves the filtrations~$C_{i,\gdot}$. This allows us to define a normal subgroup scheme $U \subset \bAut(\ul{N})$ by
$$
U(R) := \Big\{f\in \bAut(\ul{N})(R)\Bigm| \gr_a(f) = \id_{B_a}\quad \hbox{{\rm for all $a\in\cA$}}\Big\}\, ;
$$
see also [\GSAS],~5.11.

\sssection{ULemma}
{\it Lemma. ---  The identity component of~$\bAut(\ul{N})$ is contained in~$U$.}
\Dskip

\Proof By construction, $U$ is the kernel of the homomorphism
$$
h\colon \bAut(\ul{N}) \longrightarrow \dirprod_{a\in\cA} \GL(B_a)
\qquad\hbox{{\rm given by}}\quad f \mapsto \big(\gr_a(f)\big)\, .
$$
It suffices to show that~$h$ factors through a finite \'etale subgroup scheme of $\prod \GL(B_a)$.

Take~$a\in\cA$, and let~$r$ be the smallest positive integer such that $\rho^r(a) = a$. The composition
$$
\cT_a := \big(B_a \mapright{t_a} B_{\rho(a)} \mapright{t_{\rho(a)}} \;\cdots\; \mapright{t_{\rho^{r-1}(a)}} B_{\rho^r(a)} = B_a\big)
$$
is a $\Frob_k^r$-linear automorphism of~$B_a$. Then $\cB_a := \{b \in B_a\mid \cT_a(b) = b\}$ is an $\mF_{p^r}$-subspace of~$B_a$ such that the natural map $\cB_a \otimes_{\mF_{p^r}} k \to B_a$ is an isomorphism. It follows that the subgroup scheme $\Gamma \subset \prod \GL(B_a)$ given by the automorphisms that commute with all $\cT_a$ is finite \'etale. On the other hand, if $f \in \bAut(\ul{N})$ then for every~$\alpha\in\cA$ we have the relation $\gr_{\rho(\alpha)}(f) \circ t_\alpha = t_\alpha \circ \gr_\alpha(f)$. Iterating this we find that $\gr_a(f)$ commutes with~$\cT_a$. Hence $h$~factors through~$\Gamma$. \QED
\Cskip

As we shall see later, $U$ is connected; hence it is the identity component of $\bAut(\ul{N})$.

\sssection{ijj}
We describe the group scheme $U$ in more detail; see also [\GSAS],~5.11. Let $R$ be a $k$-algebra, and let $f \in \bAut(\ul{N})(R)$. Then $f$~respects the decomposition $N \otimes R = \oplus_{i\in\cI} N_i \otimes R$. Write $A_f(i,j^\prime,j)$ for the $(j^\prime,j)$-th matrix coefficient (with respect to the basis $e_{i,j}$) of the automorphism of~$N_i \otimes R$ induced by~$f$. Clearly $f$ is fully determined by the map $A_f\colon \cI \times \{1,\ldots,d\}^2 \rightarrow R$. Our task is to describe which maps $A \colon \cI \times \{1,\ldots,d\}^2 \rightarrow R$ arise as~$A_f$ for some $f\in U$. The first obvious condition is that for all $f\in U$ the map $A= A_f$ satisfies
$$
A(i,j^\prime,j) = 0\quad \hbox{{\rm if $j^\prime > j\, $,}}
\ \hbox{and}\quad 
A(i,j^\prime,j) = 1\quad \hbox{{\rm if $j^\prime = j\, $.}}
\leqno(\refn{ijj}.1)
$$ 

Consider the set $\cI \times \{1,\ldots,d\}$. Given $(i,j)$ there is a unique index~$t_i(j) \in \{1,\ldots,d\}$ such that either $F(e_{i,j}) = e_{i+1,t_i(j)}$ or $e_{i,j} = V(e_{i+1,t_i(j)})$. The~$t_i$ are permutations of $\{1,\ldots,d\}$. Given $\gf(i)$ we can easily compute $t_i$ from~$\dot w_i$ and vice versa; see the formulas for~$F$ and~$V$ in~\refn{SOnpol}. When we draw a Dieudonn\'e module as in Figure~\refn{Fig:Exa}, it is the permutations~$t_i$ that we ``see'' in the picture. Write~$T$ for the permutation of $\cI \times \{1,\ldots,d\}^2$ given by $T(i,j^\prime,j) = \big(i+1,t_i(j^\prime),t_i(j)\big)$.

Let $\ijj{0} \in \cI \times \{1,\ldots,d\}^2$. For $n \in \mZ$ write $\ijj{n} := T^n\ijj{0}$. Suppose $j_0^\prime < j_0$. We can think of the triple $\ijj{0}$ as a pair of base vectors $e_{i_0,j^\prime_0}$ and $e_{i_0,j_0}$ with $e_{i_0,j^\prime_0}$ ``below'' $e_{i_0,j_0}$. Going from $\ijj{0}$ to $\ijj{1}$ there are four possible configurations, as shown in Figure~\refn{Fig:FVdiag}. (For the interpretation of these illustrations see~\refn{picture}; cf.\ also the figures in~[\GSAS].) In cases 1 and~2 we say that $T$ is {\it parallel\/} at $\ijj{0}$; in case~3 we say $T$ is {\it down-up\/}, in case~4 we say it is {\it up-down\/}.
\midinsert
\line{\hfill\epsfxsize=.95\hsize\epsfbox{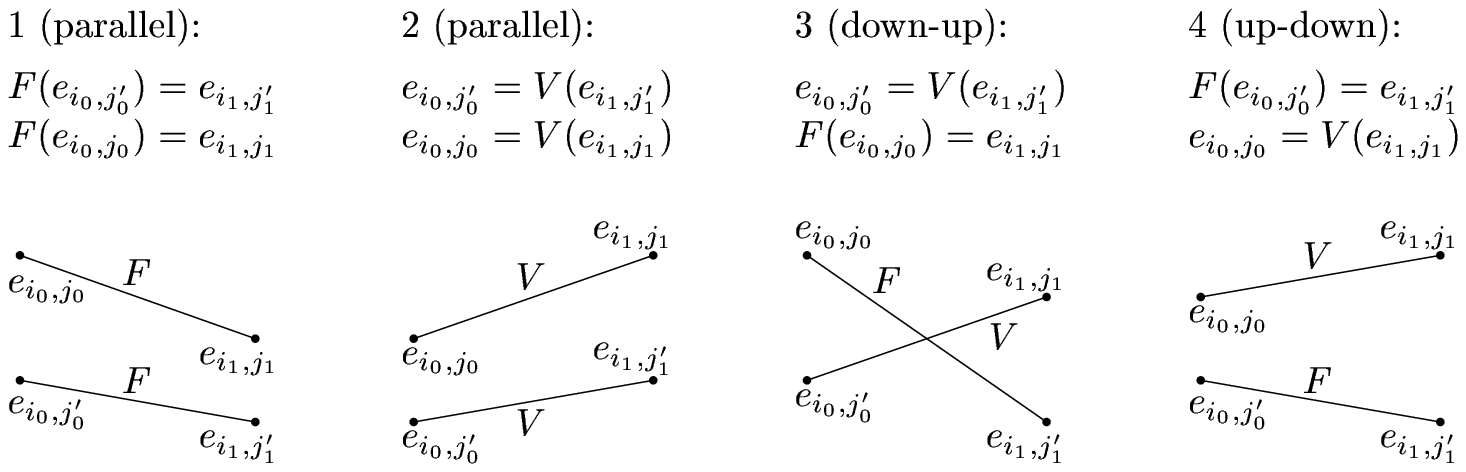}\hfill}
\figlabel{Fig:FVdiag}
\endinsert

\sssection{ALemma}
{\it Lemma. --- Suppose given a map $A \colon \cI \times \{1,\ldots,d\}^2 \rightarrow R$ that satisfies\/ {\rm (\refn{ijj}.1)}. Then $A = A_f$ for some $f \in \bAut(\ul{N})(R)$ if and only if the following two conditions are satisfied for all $\ijj{0}$ with $j^\prime_0 < j_0$:
\smallskip
\diseqitem{{\rm (\refn{ALemma}.1)}} if $T$ is parallel at $\ijj{0}$ then $A\ijj{1} = A\ijj{0}^p$;
\diseqitem{{\rm (\refn{ALemma}.2)}} if $T$ is up-down at $\ijj{0}$ then $A\ijj{0}^p = 0$.}
\Dskip

\Proof First suppose that $A = A_f$ for some $f \in \bAut(\ul{N})(R)$. That (\refn{ALemma}.1) and (\refn{ALemma}.2) hold is immediate from the fact that $f$ is compatible with the linear maps $F^\sharp$ and~$V^\flat$. (See \refn{DSpace} for notation.) For instance, suppose $T$ is up-down at $\ijj{0}$. Writing out elements on the basis $e_{i_0,1},\ldots,e_{i_0,d}$ we have 
$$
f(e_{i_0,j_0}) = e_{i_0,j_0} + \cdots + A\ijj{0} \cdot e_{i_0,j_0^\prime} + \cdots\, .
$$
Applying $F^\sharp$ to both sides, and using that $F^\sharp \circ f^{(p)} = f \circ F^\sharp$, we find a relation $0 = 0 + \cdots + A\ijj{0}^p \cdot e_{i_1,j_1^\prime} + \cdots\; $. So indeed $A\ijj{0}^p = 0$. 

Conversely, suppose $f$ is a $\kappa \otimes_{\Fp} R$-linear automorphism of~$N$ such that the corresponding matrices $A(i,j^\prime,j)$ satisfy (\refn{ijj}.1), (\refn{ALemma}.1) and (\refn{ALemma}.2). We have to verify that $f$ is compatible with $F^\sharp$ and~$V^\flat$. To check the compatibility with~$V^\flat$, let us consider a base vector $e_{i_1,j_1}$. Note that for each~$i$ the image of~$F_{i-1} \colon N_{i-1} \to N_i$ is spanned by the {\it first\/} $d-\gf(i-1)$ base vectors $e_{i,j}$. Hence if $V^\flat(e_{i_1,j_1}) = 0$ then also $V^\flat(e_{i_1,j_1^\prime}) = 0$ for all $j_1^\prime < j_1$, and it follows that $f^{(p)} V^\flat(e_{i_1,j_1}) = 0 = V^\flat f(e_{i_1,j_1})$. The other possibility is that $V^\flat(e_{i_1,j_1}) = e_{i_0,j_0}$. In this case we have
$$
\eqalign{fV^\flat(e_{i_1,j_1}) &= e_{i_0,j_0} + \sum_{j_0^\prime < j_0} A\ijj{0}^p \cdot e_{i_0,j_0^\prime}\cr
V^\flat f(e_{i_1,j_1}) &= e_{i_0,j_0} + \sum_{j_1^\prime < j_1} A\ijj{1} \cdot V(e_{i_1,j_1^\prime})\, .\cr}
$$
If $j_1^\prime < j_1$, write $j_1^\prime = t_{i_0}(j_0^\prime)$. Either $j_0^\prime > j_0$, in which case $T$ is down-up at $\ijj{0}$ and $V(e_{i_1,j_1^\prime}) = 0$, or $j_0^\prime < j_0$, in which case $T$ is parallel at $\ijj{0}$ and (\refn{ALemma}.1) gives $A\ijj{1} = A\ijj{0}^p$. Combining this we find that $f^{(p)} V^\flat = V^\flat f$. We leave it to the reader to verify, in a similar manner, that $F^\sharp \circ f^{(p)} = f \circ F^\sharp$. \QED   

\sssection{indefpar}
As before, suppose $j^\prime_0 < j_0$. It can be shown (see [\GSAS], 4.16 and 4.19) that $T$ is parallel at $\ijj{n}$ for all $n \in \mZ$ if and only $(i_0,j^\prime_0)$ and $(i_0,j_0)$ belong to the same canonical block, by which we mean that there is an index~$r$ such that $e_{i_0,j^\prime_0}$ and $e_{i_0,j_0}$ both lie in $C_{i_0,r}\backslash C_{i_0,r-1}$. Hence in order to have $A=A_f$ for some $f\in U$, the map $A$ should satisfy the additional requirement
\smallskip
\diseqitem{(\refn{indefpar}.1)} if $j^\prime_0 < j_0$ and $T$ is parallel at $\ijj{n}$ for all~$n$ then $A(i_0,j^\prime_0,j_0) = 0$.

\sssection{tracks}
To summarize our conclusions, let us define a {\it track\/} in $\cI \times \{1,\ldots,d\}^2$ to be a sequence
$$
\ijj{0},\ijj{1},\ldots,\ijj{b}
$$
with $\ijj{n} = T^n\ijj{0}$ for all~$n$, such that

\noindent
--- $j_n^\prime < j_n$ for all $0 \leq n \leq b$;

\noindent
--- $T$ is parallel at $\ijj{n}$ for all $0 \leq n \leq b-1$;

\noindent
--- $T$ is not parallel at $\ijj{b}$;

\noindent
--- $T$ is not parallel at $\ijj{{-1}} := T^{-1}\ijj{0}$.

\noindent
We call $\ijj{0}$ the start of the track, $\ijj{b}$ its end, and $b$ its length. If $T$ is up-down at $\ijj{b}$ then we call the track a {\it ud-track}; if $T$ is down-up at $\ijj{b}$ then we call it a {\it du-track}.

Suppose there are $\mu$ different du-tracks in $\cI \times \{1,\ldots,d\}^2$ and that there are $\nu$ different ud-tracks, of lengths $b_1, \ldots,b_\nu$. Then we find an isomorphism of $k$-schemes 
$$
U \isomarrow \Spec\Big(k[x_1,\ldots,x_\mu,y_1,\ldots,y_\nu]/(y_1^{p^{b_1}},\ldots,y_\nu^{p^{b_\nu}})\Big)
$$
by associating to $f \in U$ the matrix coefficients $A\ijj{0}$ for all starting points $\ijj{0}$ of a track. In particular, $U$ is connected, and together with Lemma~\refn{ULemma} it follows that $U = \bAut(\ul{N})^0$.

\sssection{PfiEnd}
To compute $\mu = \dim(U)$ we count the du-tracks by their end points. These end points are simply all triples $\ijj{{}}$ at which $T$ is down-up. Hence if for $i \in \cI$ we define
$$
\eqalign{
S(i) :=\ &\{(j^\prime,j) \mid 1 \leq j^\prime < j \leq d, F(e_{i,j^\prime}) = 0\ \hbox{and}\ F(e_{i,j}) \neq 0\}\cr
=\ &\{(j^\prime,j) \mid 1 \leq j^\prime < j \leq d\ \hbox{and}\ 
\dot w_i(j^\prime) \leq \gf(i) < \dot w_i(j)\}}
$$
then the result is that $\dim(U) = \sum_{i\in\cI} \# S(i)$. It follows from (\refn{l(w)npol}.1) that if $\dot w_i(j^\prime) = n \leq \gf(i)$ then there are $j^\prime - n$ elements in the range $1,\ldots,j^\prime$ for which $\dot w_i$ takes a value $> \gf(i)$. So there are $\big(d-\gf(i)\big) - (j^\prime - n)$ elements~$j$ such that $(j^\prime,j) \in S(i)$. This gives, using~(\refn{l(w)npol}.2):
$$
\eqalign{
\# S(i) 
&= \sum_{n=1}^{\gf(i)} \big(d-\gf(i)\big) - (\dot w_i^{-1}(n) - n)\cr
&= \gf(i) \cdot \big(d-\gf(i)\big) -  \sum_{n=1}^{\gf(i)} (\dot w_i^{-1}(n) - n) = \dim(\mX_i) - \ell(\dot w_i)\, ,}
$$
which is the formula we want. \QED\enspace (Thm.~\refn{AutDimThm}, (\romno1).)

\sssection{AutDimExa}
{\it Example.\/} --- We consider the Dieudonn\'e module illustrated in Figure~\refn{Fig:AutExa}; here $m=4$ and $d=4$, with $\gf(i)=1$, $\gf(i+1)=0$, and $\gf(i+2)=\gf(i+3) =2$.
\midinsert
\centerline{\epsfbox{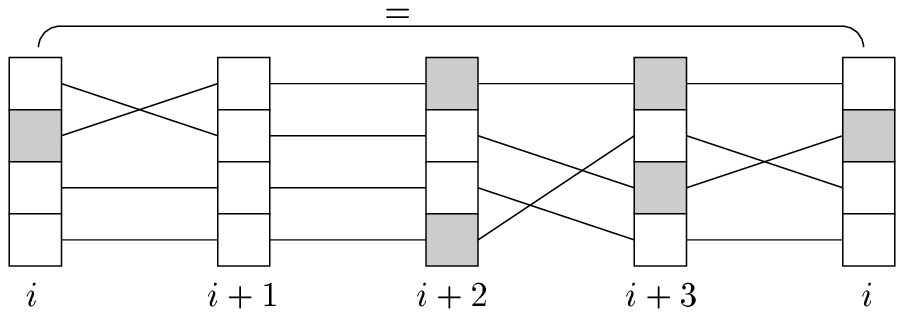}}
\figlabel{Fig:AutExa}
\endinsert
\noindent
The distinguished representative of the element~$w$ is given by
$$
\dot w_{i} = \left[\matrix{1&2&3&4\cr 2&3&1&4\cr}\right]\, ,
\dot w_{i+1} = \left[\matrix{1&2&3&4\cr 1&2&3&4\cr}\right]\, ,
\dot w_{i+2} = \left[\matrix{1&2&3&4\cr 1&3&4&2\cr}\right]\, ,
\dot w_{i+3} = \left[\matrix{1&2&3&4\cr 3&1&4&2\cr}\right]\, ,
$$
which gives $\ell(\dot w) = 2+0+2+3=7$. As $\dim(\mX) = 1\cdot 3 + 0\cdot 4 + 2\cdot 2 + 2\cdot 2 = 11$ we should find a $4$-dimensional automorphism group. Note that there are indeed precisely four du-tracks; they are the sequences
$$
\eqalignno{
&(i,1,4) \to (i+1,1,3) \to (i+2,1,3)\, ;\quad (i+3,2,4) \to (i,3,4)\, ;\quad (i+3,2,3)\cr
\noalign{\hbox{and}}
&(i+3,1,3) \to (i,1,2) \to (i+1,1,2) \to (i+2,1,2)\, .\cr}
$$
A straightforward but tedious calculation shows that an automorphism of the Dieudonn\'e module (over a field) is given by the four-tuple of matrices
$$
\left[\pmatrix{\vbox{%
\halign{\hss$#$\hss\hskip6pt&\hss$#$\hss\hskip6pt&\hss$#$\hss\hskip6pt&\hss$#$\hss\cr%
\zeta&w^p&&t\cr%
&\zeta^{p^4}&&\cr%
&&\eta&u^p\cr%
&&&\eta^{p^4}\cr}}}\, ,
\pmatrix{\vbox{%
\halign{\hss$#$\hss\hskip6pt&\hss$#$\hss\hskip6pt&\hss$#$\hss\hskip6pt&\hss$#$\hss\cr%
\zeta^p&w^{p^2}&t^p&\cr
&\zeta^{p^5}&&\cr
&&\eta^{p^5}&\cr
&&&\eta^p\cr}}}\, ,
\pmatrix{\vbox{%
\halign{\hss$#$\hss\hskip6pt&\hss$#$\hss\hskip6pt&\hss$#$\hss\hskip6pt&\hss$#$\hss\cr%
\zeta^{p^2}&w^{p^3}&t^{p^2}&\cr
&\zeta^{p^6}&&\cr
&&\eta^{p^6}&\cr
&&&\eta^{p^2}\cr}}}\, ,
\pmatrix{\vbox{%
\halign{\hss$#$\hss\hskip6pt&\hss$#$\hss\hskip6pt&\hss$#$\hss\hskip6pt&\hss$#$\hss\cr%
\zeta^{p^7}&&w&\cr
&\eta^{p^7}&v&u\cr
&&\zeta^{p^3}&\cr
&&&\eta^{p^3}\cr}}}\right]\, ,
$$
where $t$, $u$, $v$ and $w$ are chosen arbitrarily, and $\zeta$ and $\eta$ satisfy $\zeta^{p^8}=\zeta$ and $\eta^{p^8}=\eta$. So we indeed find a $4$-dimensional automorphism group. The identity component $U \subset \bAut(\ul{N})$ corresponds to setting $\zeta=\eta=1$. Note that the variables $t$, $u$, $v$ and~$w$, and the places where they occur in the above four-tuple of matrices, correspond exactly to the four du-tracks. Note further that this only describes the points of $\bAut(\ul{N})$, which is non-reduced, with values in a field; as there are seven ud-tracks the Lie algebra of $\bAut(\ul{N})$ has dimension $4+7 = 11$.

\sssection{LieAut}
Retaining the notations used above, we find that $\bAut(\ul{N})_\red^0 \cong \mA^\mu_k$ as schemes. The isomorphism is obtained by sending a point of $U = \bAut(\ul{N})_\red^0$ to the matrix coefficients $A\ijj{0}$ for all triples $\ijj{0}$ that form the start of a du-track. If $\ijj{n} = T^n\ijj{0}$ then we have seen that $A\ijj{n} = A\ijj{0}^{p^n}$. As a result, when we identify the Lie algebra of $\GL_{\kappa \otimes k}(N) = \prod_{i \in\cI} \GL_{d,k}$ with $\prod_{i\in\cI} M_d(k)$, we have
$$
\Lie\big(\bAut(\ul{N})_\red\big) = \left\{(X_i) \in \dirprod_{i\in\cI} M_{d}(k) \biggm| 
\vcenter{
\setbox0=\hbox{{\rm $X_i(j^\prime,j) \neq 0$ only if $\ijj{{}}$}}
\copy0
\hbox to\wd0{{\rm \hfill is the start of a du-track\hfill}}
}\right\}\, .
$$

\ssection{The polarized case}{Autpol}
\Bskip

\noindent
Throughout this section we work over an algebraically closed field $k$ with $\charact(k) = p >2$. Let $\cD = (B,\ast,\epsilon)$ be a triple as in section~\refn{BT1EP}. Let~$\ul{Y}$ be a BT$_1$ with $\cD$-structure over~$k$, of type $(d,\gf,\delta)$. We write $\ul{N} = (N,F,V;\phi,\iota)$ for the associated Dieudonn\'e space with $\cD$-structure and $\ul{N}^\prime := (N,F,V,\iota)$ for the Dieudonn\'e space with $B$-structure obtained by forgetting the polarization form~$\phi$. Define $(G,\mX)$ as in~\refn{pairs2}, and let $(G^0,\mX^0)$ be as in~\refn{delta}. Let $w := w(\ul{Y}) \in W_{\mX}\backslash W_G$ and $w^0 := w^0(\ul{Y}) \in W_{\mX^0}\backslash W_{G^0}$. As explained in~\refn{delta} and~\refn{ww0}, there is only a difference between $w$ and~$w^0$ if there are factors of type~D, and the two determine each other once we fix the invariant~$\delta$.

As a general notational convention, we use a prime ${}^\prime$ for objects obtained by ``forgetting the polarization form~$\phi$''. E.g., we write $G^\prime := \GL_{B \otimes_{\Fp} k}(N)$. In particular, our convention means that the objects studied in section~\refn{Autnonpol} will now appear on stage equipped with a prime. 

\sssection{AutRedLem}
{\it Lemma. --- We have $\bAut(\ul{N})_\red = \bAut(\ul{N}^\prime)_\red \cap G$.}
\Dskip

\Proof~Recall that $G = \UU_{B \otimes_{\Fp} k}(N,\phi)$. Let $A := \bAut(\ul{N}^\prime)_\red \cap G$. It is clear that $\bAut(\ul{N})_\red \subset A \subset \bAut(\ul{N})$, so it suffices to show that $A$ is smooth at the origin. Let $(R_i,\gm_i)$ be artinian local $k$-algebras with $R_i/\gm_i \cong k$. Let $\pi\colon R_1 \twoheadrightarrow R_2$ be a surjective homomorphism such that $I := \Ker(\pi)$ satisfies $\gm_1 \cdot I = 0$. Let $\alpha_2 \in A(R_2)$ be an element which reduces to the identity modulo~$\gm_2$. We have to show that $\alpha_2$ can be lifted to an element $\alpha_1 \in A(R_1)$.

As $\bAut(\ul{N}^\prime)_\red$ is a smooth $k$-group scheme we can lift $\alpha_2$ to an $R_1$-valued point $\alpha_1$ of $\bAut(\ul{N}^\prime)_\red$. If $\ga$ is the Lie algebra of $\bAut(\ul{N}^\prime)_\red$ then the set of all such liftings of~$\alpha_2$ is a principal homogeneous space under $\ga \otimes_k I$. Our task is to show that there exists an element $\beta \in \ga \otimes_k I$ such that $\tilde\alpha_1 := \alpha_1 + \beta$ preserves the form~$\phi$.

As $\phi$ satisfies (\refn{DieuTh}.2), it gives rise to an involution $g \mapsto g^\ast$ on $G^\prime = \GL_{B \otimes k}(N)$. It follows from the relation (\refn{DieuTh}.1) that this involution preserves the closed subgroup scheme $\bAut(\ul{N}^\prime)$; hence it also preserves $\bAut(\ul{N}^\prime)_\red$. It follows that $\alpha_1^\ast$ is also an $R_1$-valued point of $\bAut(\ul{N}^\prime)_\red$. By construction, $\alpha_1^\ast \alpha_1$ is congruent to the identity modulo~$I$, and it follows that $\gamma := \id - \alpha_1^\ast \alpha_1$ is an element of~$\ga \otimes_k I$.

Consider the bilinear form 
$$
\Psi\colon (N \otimes_k R_1) \times (N \otimes_k R_1) \tto R_1
\qquad
\hbox{{\rm given by}}\ \Psi(n,n^\prime) = \phi(n,n^\prime) - \phi(\alpha_1 n,\alpha_1 n^\prime)\, .
$$
Because $\alpha_1$ lifts $\alpha_2$ and $\alpha_2$ preserves~$\phi$, the form $\Psi$ takes values in~$I$. Further, $\Psi(n,n^\prime)$ only depends on the classes of~$n$ and~$n^\prime$ modulo~$\gm_1$, as $\gm_1 \cdot I = 0$. In sum, $\Psi$ gives rise to a bilinear form $\psi \colon N \times N \tto I$. The element $\gamma \in \ga \otimes_k I \subset \End_{B \otimes k}(N) \otimes_k I$ is the unique element such that $\psi(n,n^\prime) = \phi(\gamma n,n^\prime)$ for all $n$, $n^\prime \in N$. Note that $\gamma^\ast = \gamma$. 

If $\tilde\alpha_1 := \alpha_1 + \beta$ for some $\beta \in \End_{B \otimes k}(N) \otimes_k I$ then the corresponding forms $\psi$ and $\tilde\psi$ are related by
$$
\tilde\psi(n,n^\prime) = \psi(n,n^\prime) - \phi\big((\beta^\ast \bar\alpha_1 +\bar\alpha_1^\ast \beta)n,n^\prime\big)\, ,
$$
where $\bar\alpha_1$ is the reduction of $\alpha_1$ modulo~$\gm_1$. By assumption we have $\bar\alpha_1 = \id_N$. As $\gamma^\ast = \gamma$ and $\charact(k) \neq 2$ we can take $\beta = \gamma/2$. This gives a lifting of~$\alpha_2$ as desired. \QED

\sssection{BasCas}
To prove (\romno2) of Thm.~\refn{AutDimThm} it suffices to consider the following four cases. In each case we further describe the possibilities for the type $(d,\gf,\delta)$, with notation as in section~\refn{BT1EP}.
\smallskip

\typesitem{Case C:} $B = \tilde\kappa =\kappa$ is a finite field with $\ast=\id$ and $\epsilon = -1$. The type $(d,\gf)$ is fully determined by a number $q\in \mZ_{> 0}$ via $d=2q$ and $\gf(i)=q$ for all $i\in\cI = \Hom(\kappa,k)$.

\typesitem{Case D:} $B = \tilde\kappa =\kappa$ is a finite field with $\ast=\id$ and $\epsilon = +1$. The type $(d,\gf)$ is as in Case~C. Further we have an invariant $\delta$, which  is an arbitrary function $\cI \to \mZ/2\mZ$.

\typesitem{Case AU:} $B=\tilde\kappa$ is a finite field, with $\ast$ an automorphism of order~$2$ and $\epsilon =1$. 
Let $\kappa$ be the fixed field of~$\ast$. The type $(d,\gf)$ consists of a positive integer~$d$ and a function $\gf\colon\cItil \to \mZ_{\geq 0}$ with $\gf(\tau) + \gf(\bar\tau) = d$ for all $\tau\in\cItil$.

\typesitem{Case AL:} $B = \tilde\kappa = \kappa \times \kappa$ with $\kappa$ a finite field, with $(x_1,x_2)^\ast = 
(x_2,x_1)$, and $\epsilon = 1$. The type $(d,\gf)$ is as in Case~AU. Note that in this case $\cItil = \cI \amalg \cI$, where the two copies of~$\cI$ are interchanged under $\tau \mapsto \bar\tau$; hence $\gf$ is determined by its restriction to the first copy of~$\cI$, which is a function $\gf\colon \cI \to \mZ_{\geq 0}$ with $\gf(i) \leq d$ for all $i\in\cI$.

\noindent
The reduction to these basic cases is an application of Morita equivalence. See also [\GSAS],~5.2 and~6.2. Note that the labels introduced in~\refn{CDAfact} correspond, via Morita equivalence, to the above labels, with type~A subdivided into two cases.

In Case~AL, every BT$_1$ with $(B,\ast,\epsilon)$-structure is of the form $\ul{Y} \times \ul{Y}^D$ where $\ul{Y}$ is a BT$_1$ with $\kappa$-structure, and where the $\epsilon$-duality $\lambda\colon \ul{Y} \times \ul{Y}^D \to \ul{Y}^D \times \ul{Y}$ is given by $(y_1,y_2) \mapsto (y_2,\epsilon y_1)$. This reduces Case~AL to the study of non-polarized BT$_1$, and one checks without difficulty that (\romno2) of the theorem follows in this case from~(\romno1).

\sssection{PfCD}
Suppose we are in Case~C or Case~D. We have $G = \prod_{i\in\cI} G_i$ with $G_i \cong \Sp_{2q,k}$ (Case~C) or $G_i \cong \OO_{2q,k}$ (Case~D). The Weyl group of $G_i$ can be described as
$$
W_{G_i} = \mH_q := \Big\{\sigma \in \gS_{2q} \Bigm| \sigma(j) + \sigma(2q+1-j) = 2q+1\quad \hbox{{\rm for all $j$}}\Big\}\, .\leqno(\refn{PfCD}.1)
$$
In each factor $W_{G_i}$ we take as generators the $q-1$ elements of the form $(j\; j+1)(2g+1-j\; 2g-j)$ together with the transposition $(q\; q+1)$. The subgroup $W_\mX \subset W_G$ is the product of the groups
$$
W_{\mX_i} = \Big\{\sigma\in W_{G_i} \Bigm| \sigma\{1,\ldots,q\} = \{1,\ldots,q\}\Big\} \cong \gS_q\, .
$$

Consider the element $w \in W_\mX\backslash W_G$. Write $\dot w = (\dot w_i)_{i\in\cI}$ for its distinguished representative. The length of $\dot w$ is given by $\ell(\dot w) = \sum_{i\in\cI} \ell(\dot w_i)$, with
$$
\ell(\dot w_i) = \sum_{n=1}^{l_i} \big(q+1-\dot w_i^{-1}(q+n)\big)\, ;\leqno(\refn{PfCD}.2)
$$
here $l_i \in \{1,\ldots,q\}$ is the largest index for which $\dot w_i^{-1}(q+l_i) \leq q$.

In Case~D the structure of the Dieudonn\'e module $\ul{N}$ is easiest to describe in terms of the classifying element $w\in W_\mX\backslash W_G$. However, in the result that we want to prove it is not the length of~$\dot w$ that matters but rather the length of the distinguished representative of $w^0 := w^0(\ul{Y}) \in W_{\mX^0}\backslash W_{G^0}$. The length of $\dot w^0$ can be expressed directly in terms of~$\dot w$; the result is that
$$
\ell(\dot w^0) = \sum_{i\in\cI} \ell(\dot w^0_i)
\quad\hbox{with}\quad
\ell(\dot w^0_i) = \sum_{n=1}^{l_i} \big(q-\dot w_i^{-1}(q+n)\big)\, ;\leqno(\refn{PfCD}.3)
$$
here $l_i \in \{1,\ldots,q\}$ is again the largest index for which $\dot w_i^{-1}(q+l_i) \leq q$. Let us further notice that, still in Case~D, the invariant $\delta\colon \cI \to \mZ/2\mZ$ can be read from the element $\dot w$ by the formula
$$
\delta(i) = \#\big\{j \leq q \bigm| \dot w_i(j) > q\big\} \bmod 2\, .
$$

Similar to what was done in~\refn{SOnpol} we have an explicit description of~$\ul{N}$ in terms of the element~$\dot w$. Namely, we can choose a $k$-basis $e_{i,j}$ for~$N$, for $i\in\cI$ and $j \in \{1,\ldots,2q\}$, such that $a \in \kappa$ acting on $e_{i,j}$ as multiplication by $i(a) \in k$, and such that $F$ and $V$ are given by
$$
\eqalign{
F_i(e_{i,j}) = e_{i+1,n}
\quad\hbox{and}\quad
V(e_{i+1,n}) = 0\qquad &\hbox{if $\dot w_i(j) = q+n$;}\cr
F_i(e_{i,j}) = 0
\quad\hbox{and}\quad
V(e_{i+1,n}) = e_{i,j}\qquad &\hbox{if $\dot w_i(j) = -q+n$.}\cr}
$$
The $\epsilon$-symmetric form $\phi$ on~$N$ is the orthogonal sum of forms $\phi_i \colon N_i \times N_i \to k$. With respect to the basis $\{e_{i,j}\}_{j=1,\ldots,2q}$ the form $\phi_i$ is given by an (invertible) anti-diagonal matrix~$\Phi_i$. In Case~D we can choose our basis such that $\Phi_i = \antidiag(1,\ldots,1)$. In Case~C there is not, in general, a natural choice of a normal form for the~$\Phi_i$. Note that (\refn{DieuTh}.1) gives a relation between $\Phi_i$ and~$\Phi_{i+1}$. Any collection of anti-diagonal matrices $\{\Phi_i\}_{i\in\cI}$ satisfying the relations imposed by (\refn{DieuTh}.1) gives a Dieudonn\'e module with $(\kappa,\id,-1)$-structure, and up to isomorphism this object is independent of the chosen collection~$\{\Phi_i\}$. We refer to [\GSAS],~5.8 for details.

Via the chosen ordered bases $\{e_{ij}\}$ we may identify the Lie algebra of $G^\prime = \prod_{i\in\cI} \GL_{2q,k}$ with $\prod_{i\in\cI} M_{2q}(k)$. As we have seen in~\refn{LieAut}, the Lie algebra of~$\bAut(\ul{N}^\prime)_\red$ is given by
$$
\Lie\big(\bAut(\ul{N}^\prime)_\red\big) = \left\{(X_i) \in \dirprod_{i\in\cI} M_{2q}(k) \biggm| 
\vcenter{
\setbox0=\hbox{{\rm $X_i(j^\prime,j) \neq 0$ only if $\ijj{{}}$}}
\copy0
\hbox to\wd0{{\rm \hfill is the start of a du-track\hfill}}
}\right\}\, .
$$
Note that if we subdivide each matrix $X_i$ into four blocks of size $q \times q$, say
$$
X_i = \pmatrix{\alpha_i&\beta_i\cr \gamma_i & \delta_i}\, ,
$$
then for $(X_i)_{i \in \cI} \in \Lie\big(\bAut(\ul{N}^\prime)_\red\big)$, only the blocks $\beta_i$ can be nonzero; this corresponds to the fact that if $\ijj{{}}$ is the start of a track then $j^\prime \leq \gf(i) = q < j$.

Write $J = J_q$ for the anti-diagonal matrix of size $q\times q$ with all anti-diagonal coefficients equal to~$1$. Given a square matrix~$A$ of size $q\times q$, write ${}^s\! A = J\cdot {}^t\! A \cdot J$ for its reflection in the anti-diagonal. So $\big({}^s\! A\big)_{t,u} = A_{q+1-u,q+1-t}$.

In the symplectic case (Case~C) there are diagonal matrices~$C_i$ such that
$$
\Phi_i = \pmatrix{C_i&0\cr 0&1} \cdot 
\pmatrix{0&J\cr-J&0} \cdot
\pmatrix{C_i&0\cr 0&1}\, .
$$
(Here we write matrices in block form.) One computes that the Lie algebra  of the symplectic group $\Sp(N_i,\phi_i)$ is then given by
$$
\gsp(N_i,\phi_i) = \left\{\pmatrix{\alpha&\beta\cr \gamma&\delta} \in M_{2q}(k) \biggm|
\vcenter{
\setbox0=\hbox{$C_i \cdot \alpha\cdot C_i^{-1} = -{}^s\!\delta\, ,\ \beta = C_i^{-1}\cdot {}^s\beta \cdot {}^s C_i$}
\copy0
\hbox to\wd0{\hfil and ${}^s\gamma = {}^s C_i \cdot \gamma\cdot C_i^{-1}$\hfil}}
\right\}\, .
$$

In the orthogonal case (Case~D), $\Phi_i = \antidiag(1,\ldots,1)$. This gives
$$
\go(N_i,\phi_i) = \left\{\pmatrix{\alpha&\beta\cr \gamma&\delta} \in M_{2q}(k) \biggm|
\vcenter{
\setbox0=\hbox{${}^s\!\beta = -\beta\, , {}^s\!\gamma = -\gamma$}
\copy0
\hbox to\wd0{\hfil and $\alpha + {}^s\!\delta = 0$\hfil}}
\right\}\, .
$$ 

Given an element $\ijj{{}} \in \cI \times \{1,\ldots,d\}^2$, write ${}^s\!\ijj{{}} := (i,2q+1-j,2q+1-j^\prime)$. Observe that ${}^s\!\big({}^s\!\ijj{{}}\big) = \ijj{{}}$. If
$$
\gt:\quad \ijj{0},\ijj{1},\ldots,\ijj{b}
$$
is a track, write
$$
{}^s\gt :=\quad {}^s\!\ijj{0},{}^s\!\ijj{1},\ldots,{}^s\!\ijj{b}\, .
$$
We claim that ${}^s\gt$ is again a track in $\cI \times \{1,\ldots,d\}^2$. This is a direct consequence of the definitions once we know that the permutation~$t_i$ defined in~\refn{ijj} is an element of the group $\mH_q$ defined in (\refn{PfCD}.1). Using that $\dot w_i \in \mH_q$, this can be checked using the formulas for $F$ and~$V$ given above.

We can now complete the proof of \refn{AutDimThm}, (\romno2) in Cases~C and~D. First we do Case~C. By Lemma~\refn{AutRedLem} the dimension of $\bAut(\ul{N})$ equals the dimension of $\Lie\big(\bAut(\ul{N}^\prime)_\red\big) \cap \gsp_{\kappa \otimes k}(N,\phi)$. By what was explained above this dimension equals the number of tracks modulo the equivalence relation $\gt \sim {}^s\gt$. (The equations for $\gsp(N_i,\phi_i)$ give one relation for every pair $\{\gt,{}^s \gt\}$ with $\gt \neq {}^s \gt$.)

As in \refn{PfiEnd} we count the tracks by their end points. These are the triples $\ijj{{}}$ with $1 \leq j^\prime < j \leq 2q$ and $\dot w_i(j^\prime) \leq q < \dot w_i(j)$. Precisely one of the tracks $\gt$ and ${}^s\gt$ has an end point $\ijj{{}}$ with $j \leq 2q+1-j^\prime$, so we impose this as an extra condition. Writing
$$
S_{{\rm C}}(i) := \{(j^\prime,j) \mid 1 \leq j^\prime < j \leq 2q+1-j^\prime\ \hbox{and}\ \dot w_i(j^\prime) \leq q < \dot w_i(j)\}
$$
this gives
$$
\dim\big(\bAut(\ul{N})\big) = \sum_{i \in \cI} \# S_{{\rm C}}(i)\, .
$$

Suppose $j^\prime \leq q$ is an index with $\dot w_i(j^\prime) = n \leq q$. We subdivide the set $\{1,\ldots,2q\}$, and in each range we count the numbers of times $\dot w_i$ takes a value $> q$:
$$
\vbox{\offinterlineskip
\hrule
\halign{&\vrule#&\strut\quad\hfil$#$\hfil\quad\cr
height2pt&\omit&&\omit&&\omit&\cr
& \hbox{range} && \hbox{$\# r$ with $\dot w_i(r) \leq q$} && \hbox{$\# r$ with $\dot w_i(r) > q$} & \cr
height2pt&\omit&&\omit&&\omit&\cr
\noalign{\hrule}
height2pt&\omit&&\omit&&\omit&\cr
& 1,\ldots,j^\prime-1 && n-1 && j^\prime-n &\cr
height2pt&\omit&&\omit&&\omit&\cr
\noalign{\hrule}
height2pt&\omit&&\omit&&\omit&\cr
& j^\prime && 1 && 0 &\cr
height2pt&\omit&&\omit&&\omit&\cr
\noalign{\hrule}
height2pt&\omit&&\omit&&\omit&\cr
& j^\prime+1,\ldots,2q-j^\prime && q-j^\prime && q-j^\prime &\cr
height2pt&\omit&&\omit&&\omit&\cr
\noalign{\hrule}
height2pt&\omit&&\omit&&\omit&\cr
& 2q+1-j^\prime && 0 && 1 & \cr
height2pt&\omit&&\omit&&\omit&\cr
\noalign{\hrule}
height2pt&\omit&&\omit&&\omit&\cr
& 2q+2-j^\prime,\ldots,2q && j^\prime-n && n-1 & \cr
height2pt&\omit&&\omit&&\omit&\cr}
\hrule}
$$
Thus, we find that there are $q-j^\prime$ indices~$j$ for which $(j^\prime,j) \in S_{{\rm C}}(i)$. If, as above, we let $l_i \in \{1,\ldots,q\}$ be the largest number for which $\dot w_i^{-1}(q+l_i) \leq q$, we get
$$
\# S_{{\rm C}}(i) = \sum_{n=1}^{q-l} q+1 - \dot w_i^{-1}(n)\, .
$$
Combined with (\refn{PfCD}.2) this gives
$$
\eqalign{
\# S_{{\rm C}}(i) + \ell(\dot w_i) &=
\sum_{n=1}^{q-l} \big(q+1 - \dot w_i^{-1}(n)\big) + \sum_{n=1}^l \big(q+1-\dot w_i^{-1}(q+n)\big)\cr
&= q(q+1) - \sum_{\mu = 1}^q \mu \;=\; q(q+1)/2 \;=\; \dim(\mX_i)\, .\cr}
$$
This proves our claim that $\dim\big(\bAut(\ul{N})\big) = \dim(\mX) - \ell(\dot w)$. 

In the orthogonal case the computation is similar. This time we have to count only the pairs of tracks $\{\gt,{}^s\gt\}$ with $\gt \neq {}^s \gt$. Setting
$$
S_{{\rm D}}(i) := \{(j^\prime,j) \mid 1 \leq j^\prime < j < 2q+1-j^\prime\ \hbox{and}\ \dot w_i(j^\prime) \leq q < \dot w_i(j)\}
$$
we get that $\dim\big(\bAut(\ul{N})\big) = \sum_{i \in \cI} \# S_{{\rm D}}(i)$. Counting as above and using (\refn{PfCD}.3) gives
$$
\eqalign{
\# S_{{\rm D}}(i) + \ell(\dot w^0_i) &= \sum_{n=1}^{l_i} q - \dot w_i^{-1}(q+n) + \sum_{n=1}^{q-l_i} q-\dot w_i^{-1}(n)\cr
&= q^2 - \sum_{\mu =1}^q \mu = q(q-1)/2 = \dim(\mX_i)\, ,\cr}
$$
which is what we want. \QED\enspace (Thm.~\refn{AutDimThm}, (\romno2), Cases~C and~D.)

\sssection{PfA}
Finally we consider Case~AU. Recall that we write $\cItil$ for the set of embeddings of $\tilde\kappa$ into~$k$. We have a natural $2:1$ map $\res\colon \cItil \to \cI$. For $\tau \in \cItil$, let $\bar\tau := \tau \circ \ast$. The type $(d,\gf)$ is given by an integer~$d$ and a function $\gf\colon\cItil \to \mZ_{\geq 0}$ with the property that $\gf(\tau) + \gf(\bar\tau) = d$ for all~$\tau$.

Given $(d,\gf)$, the corresponding pair of Weyl groups $W_\mX \subset W_G$ can be described as follows. Let $w_0 \in \gS_d$ be the permutation of order~$2$ given by $w_0(j) = d+1-j$. For $\pi \in \gS_d$, write $\check\pi := w_0\circ \pi \circ w_0$. Then
$$
W_G = \Big\{(\pi_\tau) \in \dirprod_{\tau\in\cItil} \gS_d \Bigm| \check\pi_\tau = \pi_{\bar\tau}\quad \hbox{{\rm for all~$\tau$}}\Big\}\, ,
$$
and
$$
W_\mX = \big\{(\pi_\tau) \in W_G \bigm| \pi_\tau\{1,\ldots,\gf(\tau)\} = \{1,\ldots,\gf(\tau)\}\quad \hbox{{\rm for all~$\tau$}}\big\}\, .
$$

Let $\dot w = (\dot w_\tau)_{\tau\in\cItil}$ for the distinguished representative of the coset~$w = w(\ul{Y})$. Let $\cR \subset \cItil$ be a subset such that precisely one member of each pair $\{\tau,\bar\tau\}$ is in~$\cR$. For each pair $\{\tau,\bar\tau\}$ we have 
$$
\ell(\dot w_\tau) = \ell(\dot w_{\bar\tau}) = \sum_{n=1}^{\gf(\tau)} \big(\dot w_\tau^{-1}(n) - n\big) = \sum_{n=\gf(\tau)+1}^d \big(n - \dot w_\tau^{-1}(n)\big)\, ,\leqno(\refn{PfA}.1)
$$
and $\ell(\dot w) = \sum_{\tau\in\cR} \ell(\dot w_\tau)$. The dimension of~$\mX$ is given by
$$
\dim(\mX) = \sum_{\tau\in\cR} \gf(\tau) \cdot \big(d-\gf(\tau)\big)\, .\leqno(\refn{PfA}.2)
$$
(Observe that $\gf(\tau) \cdot \big(d-\gf(\tau)\big) = \gf(\bar\tau) \cdot \big(d-\gf(\bar\tau)\big)$ because $\gf(\tau) + \gf(\bar\tau) = d$.)

The description of the Dieudonn\'e module $\ul{N}$ is much the same as in the previous cases. Namely, we can choose a $k$-basis $e_{\tau,j}$ for $\tau\in\cItil$ and $j\in\{1,\ldots,d\}$ such that $a\in \tilde\kappa$ act on~$e_{\tau,j}$ as multiplication by $\tau(a) \in k$, and such that Frobenius and Verschiebung are given by
$$
\eqalign{
F_\tau(e_{\tau,j}) = e_{\tau+1,n}
\quad\hbox{and}\quad
V(e_{\tau+1,n}) = 0\qquad &\hbox{if $\dot w_\tau(j) = \gf(\tau)+n$;}\cr
F_i(e_{\tau,j}) = 0
\quad\hbox{and}\quad
V(e_{\tau+1,n}) = e_{\tau,j}\qquad &\hbox{if $\dot w_\tau(j) = -d+\gf(\tau)+n$.}\cr}
$$
The form $\phi$ decomposes as an orthogonal sum of hermitian forms $\phi_i$ on $N_\tau \times N_{\bar\tau}$ (where $i = \res(\tau) = \res(\bar\tau)$). We can choose the ordered basis $e_{\tau,1}, \ldots ,e_{\tau,d},e_{\bar\tau,1}, \ldots,e_{\bar\tau,d}$ such that $\phi_i$ is given by an invertible anti-diagonal matrix~$\Phi_i$. As in Case~C there is, in general, not a natural choice of a normal form for these matrices.

For $(\tau,j^\prime,j) \in \cItil \times \{1,\ldots,d\}^2$, write ${}^s\!(\tau,j^\prime,j) := (\bar\tau,d+1-j,d+1-j^\prime)$. Following the same procedure as in~\refn{PfCD}, if $\gt$ is a track we get a ``mirrored'' track~${}^s\gt$. The condition that the form~$\psi$ is preserved gives one relation for each pair of tracks $\{\gt,{}^s\gt\}$. We find, counting the tracks by their end points, that
$$
\dim\big(\bAut(\ul{N})\big) = \sum_{\tau\in\cR} \# S_{{\rm A}}(\tau)\, ,
$$
with 
$$
S_{{\rm A}}(\tau) = \{(j^\prime,j) \mid 1\leq j^\prime < j \leq d\ \hbox{and}\ \dot w_\tau(j^\prime) \leq \gf(\tau) < \dot w_\tau(j)\}\, . 
$$
If $\dot w_\tau(j^\prime) = n \leq \gf(\tau)$ then there are $d-\gf(\tau) + n - j^\prime$ values of~$j$ with $(j^\prime,j) \in S_{{\rm A}}(\tau)$. Hence
$$
\# S_{{\rm A}}(\tau) = \sum_{n=1}^{\gf(\tau)} \big(d-\gf(\tau) + n - \dot w_\tau^{-1}(n)\big) = \gf(\tau) \big(d-\gf(\tau)\big) + \sum_{n=1}^{\gf(\tau)} n - \dot w_\tau^{-1}(n)\, .
$$
Combining this with (\refn{PfA}.1) and (\refn{PfA}.2) we obtain the equality $\ell(\dot w) + \dim\big(\bAut(\ul{N})\big) = \dim(\mX)$ that we wanted to prove. \QED\enspace (Thm.~\refn{AutDimThm}, (\romno2).)
\section{Application to Ekedahl-Oort stratifications}{}

\ssection{The Ekedahl-Oort stratification on moduli spaces of PEL type}{}

\sssection{PELdata}
We consider a moduli problem of PEL type with good reduction at a prime $p>2$. The data involved are the following.

\item{---} $(\cB,\ast)$ is a finite dimensional semi-simple $\mQ$-algebra with a positive involution;
\item{---} $\cV$ is a finitely generated faithful left $\cB$-module;
\item{---} $\phi\colon \cV \times \cV \rightarrow \mQ$ is a symplectic form ($\mQ$-bilinear, alternating and perfect) 
with the property that $\phi(bv_1,v_2) = \phi(v_1,b^\ast v_2)$ for all $b \in \cB$ and $v_1$, $v_2 \in \cV$;
\item{---} $p$ is a prime number $>2$ such that $\cB \otimes \Qp$ is unramified, i.e., isomorphic to a product of 
matrix algebras over unramified field extensions of~$\Qp$;
\item{---} $O_\cB$ is a $\mZ_{(p)}$-order in~$\cB$, stable under~$\ast$, such that $O_\cB \otimes \Zp$ is a maximal 
order in $\cB \otimes \Qp$;
\item{---} $\Lambda \subset \cV \otimes \Qp$ is a $\Zp$-lattice which is also an $O_\cB$-submodule, such that $\phi$ 
induces a perfect pairing $\Lambda \times \Lambda \to \Zp$;
\item{---} $\cG := \CSp(\Lambda,\phi) \cap \GL_{O_\cB \otimes \Zp}(\Lambda)$ is the (not necessarily connected) 
reductive group over~$\Zp$ given by the symplectic similitudes of $(\Lambda,\phi)$ that commute with the action 
of~$O_\cB$;
\item{---} $\cX$ is a $\cG(\mR)$-conjugacy class of homomorphisms $\mS \to \cG_\mR$ (with $\mS := \Res_{\mC/\mR} 
\mG_m$) that define a Hodge structure of type $(-1,0) + (0,-1)$ on~$\cV_\mR$ for which either $2\pi i\cdot \phi$ or 
$-2\pi i\cdot \phi$ is a polarization form;
\item{---} $\gc$ is the $\cG(\mC)$-conjugacy class of cocharacters of~$\cG_\mC$ associated to~$\cX$; concretely, if 
$h\in\cX$ then we have a cocharacter $\mu=\mu_h$ through which $z\in\mC^\times$ acts on $\cV^{-1,0}$ (resp.\ 
$\cV^{0,-1}$) as multiplication by~$z$ (resp.\ by~$1$);
\item{---} $E$ is the reflex field, i.e., the field of definition of the conjugacy class~$\gc$.

\sssection{ADdef}
Fix data $\cD = (\cB,\ast,\cV,\phi,O_\cB,\Lambda,\cX)$ as in \refn{PELdata}. Let $\Qbar$ be the algebraic closure 
of~$\mQ$ inside~$\mC$. We fix an embedding $\Qbar \to \Qbar_p$. Let $v$ be the corresponding place of~$E$ above~$(p)$. 
We write $O_{E,v}$ for the localization of~$O_E$ at~$v$.

Let $C_p := \cG(\Zp)$. Let $C^p$ be a compact open subgroup of $\cG(\mA^p_f)$, and put $C := C_p \times C^p$. We 
consider the moduli problem $\cA_{\cD,C}$ over $\Spec(O_{E,v})$ defined by Kottwitz in [\Kott],~\S5. If $T$ is a 
locally noetherian $O_{E,v}$-scheme then the $T$-valued points of $\cA_{\cD,C}$ are the isomorphism classes of 
four-tuples $\ul{A} = (A,\bar\lambda,\iota,\bar\eta)$ with
\item{---} $A$ an abelian scheme up to prime-to-$p$ isogeny over~$T$;
\item{---} $\bar\lambda \in \big(\NS(A) \otimes \Zp\big)/\mZ_p^\times$ the class of a prime-to-$p$ polarization;
\item{---} $\iota \colon O_\cB \to \End_T(A) \otimes \mZ_{(p)}$ a homomorphism of $\mZ_{(p)}$-algebras with 
$\iota(b^\ast) = \iota(b)^\dagger$; here $\dagger$ is the Rosati involution associated to~$\bar\lambda$;
\item{---} $\bar\eta$ a level structure of type $C^p$ on~$A$;

\noindent
such that a certain determinant condition is satisfied. For precise details we refer to Kottwitz~[\Kott],~\S5. If $C^p$ 
is sufficiently small, which we from now on assume, then $\cA_{\cD,C}$ is representable by a smooth quasi-projective 
$O_{E,v}$-scheme.

\sssection{modpdata}
Let $B := O_\cB/pO_\cB$. The involution~$\ast$ induces an involution of~$B$ which we denote by the same symbol. Let 
$\tilde{K}$ be the centre of~$\cB$, let $K := \{z \in \tilde K \mid z^\ast =z\}$, and let $O_{\tilde{K}}$ and $O_K$ be 
the respective rings of integers. Then $O_{\tilde{K}}/pO_{\tilde{K}} = \tilde\kappa$ is the centre of~$B$ and $O_K/pO_K 
= \kappa := \{z \in \tilde\kappa \mid z^\ast = z\}$.

Let $\mQ_p^\nr \subset \Qbar_p$ be the maximal unramified extension of~$\Qp$. We write $\Fpbar$ for its residue field. 
The assumptions on the data in~\refn{PELdata} imply that $v$ is an unramified prime, so $E_v \subset \mQ_p^\nr$. 

Consider an algebraically closed $k$ containing~$\Fpbar$. As usual we write $\cItil$ (resp.~$\cI$) for the set of 
embeddings of~$\tilde\kappa$ (resp.~$\kappa$) into~$k$. Using that $v$ is unramified we get, via the chosen embedding 
$\Qbar \to \Qbar_p$, natural identifications $\cItil = \Hom(\tilde{K},\mC)$ and $\cI = \Hom(K,\mC)$. Write $\cI = \cI_1 
\cup \cdots \cup \cI_\nu$ for the partition of~$\cI$ corresponding to the decomposition of~$\kappa$ as a product of 
finite fields.

The choice of $h \in \cX$ gives $\cV$ a $\mQ$-Hodge structure of type $(-1,0) + (0,-1)$ equipped with an action 
of~$\cB$. In particular we have a pair a $\cB \otimes_\mQ \mC$-modules $\cV_\mC^{-1,0} \subset \cV_\mC$. Via the above 
identifications this gives us a type~$(d,\gf)$ as in~\refn{pairs2} with $d$ constant on each of the subsets $\cI_n 
\subset \cI$ and such that (\refn{pairs2}.1) holds. This type is independent of the chosen~$h$.

A $k$-valued point of $\cA_{\cD,C}$ gives rise to a BT$_1$ with $(B,\ast,-1)$-structure $\ul{Y}$ by taking $Y := A[p]$ 
equipped with its $B$-action, and with $\lambda \colon Y \isomarrow Y^D$ the $-1$-duality induced by the 
prime-to-$p$-polarization on~$A$. The determinant condition mentioned in~\refn{ADdef} implies that $\ul{Y}$ has type 
equal to~$(d,\gf)$.

\sssection{deltalc}
{\it Lemma. --- Let $S$ be a scheme of characteristic $p > 2$. Let $\ul{Y}$ be a BT$_1$ with 
$(B,\ast,\epsilon)$-structure over~$S$. Let $k$ be an algebraically closed field, $\charact(k) = p$. For $s \in S(k)$ 
let $\delta_s \colon \cI^{{\rm D}} \to \mZ/2\mZ$ be the invariant associated to~$\ul{Y}_s$ as in\/~{{\rm 
\refn{delta}}}. Then the function $s \mapsto \delta_s$ is locally constant on~$S(k)$.}
\Dskip

\Proof~It suffices to consider the case that $S$~is the spectrum of a perfect valuation ring~$A$. Write $Q$ for the 
fraction field of~$A$ and $\gm_A$ for its maximal ideal. By Berthelot~[\Berth], we have a Dieudonn\'e theory that 
generalizes the theory over a perfect field. (In fact, the Dieudonn\'e functor is again given by taking homomorphisms 
into the Witt covectors.) So we have a $5$-tuple $\ul{N} = (N,F,V,\iota,\phi)$ as in~\refn{DieuTh} but now with $N$ a 
free $A$-module of finite rank. After base change to~$Q$ (resp.\ to~$A/\gm_A$) we retrieve the Dieudonn\'e module of 
the generic (resp.\ special) fibre of~$\ul{Y}$. 

We know that $\Ker(F_Y) \subset Y$ and $\Ker(V_Y) \subset Y^{(p/S)}$ are free subgroup schemes; their Dieudonn\'e 
modules can be identified with $\Ker(F_N) \subset N$ and $\Ker(V_N) \subset N$, respectively. The exactness of the 
sequences
$$
1\tto \Ker(F_Y) \tto Y \tto \Ker(V_Y) \tto 1
\qquad\hbox{and}\qquad
1 \tto \Ker(V_Y) \tto Y \tto \Ker(F_Y) \tto 1
$$
then implies, using [\Berth], Prop.~2.4.1, that $\Ker(F_N)$ and $\Ker(V_N)$ are direct summands of~$N$.

Using Morita equivalence one reduces the problem to the case where $B = \kappa$ is a finite field with $\ast = 
\id_\kappa$ and $\epsilon = +1$. Then the assertion reduces to the following fact: Let $N$ be a free $A$-module of even 
rank~$2q$ equipped with a perfect symmetric pairing $\phi\colon N \times N \rightarrow A$. Let $L_1$ and $L_2$ be 
direct summands of $N$ of rank~$q$ that are totally isotropic with respect to~$\phi$. Then $\dim\big((L_1 \otimes Q) 
\cap (L_2 \otimes Q)\big)$ has the same parity as $\dim\big((L_1 \otimes A/\gm_A) \cap (L_2 \otimes A/\gm_A)\big)$. 
\QED

\sssection{EOstrat}
Write $\cA_0 := \cA_{\cD,C} \otimes \Fpbar$. The lemma implies that we can decompose~$\cA_0$ as a disjoint union of 
open and closed subscheme according to the value of~$\delta$, say $\cA_0 = \amalg\, \cA_{0,\delta}$. 

Fix~$\delta$. To the type $(d,\gf,\delta)$ we associate a pair $(G^0,\mX^0)$ as in~\refn{delta}. Note that, up to a 
central factor, $G^0$ can be identified with the special fibre of~$\cG$, and that $\mX^0$ is closely related 
to~$\cX^0$. We shall not attempt to make this more precise. Now our classification results of BT$_1$ with additional 
structure give rise to a generalized Ekedahl-Oort stratification
$$
\cA_{0,\delta} = \disunion_{w \in W_{\mX^0}\backslash W_{G^0}} \cA_{0,\delta}(w)\, .
$$
More precisely, we have a partition of~$\cA_{0,\delta}$ into a disjoint union of locally closed subspaces such that the 
Zariski closure of each stratum is a union of strata. These properties were proven by Wedhorn in~[\Wedh], generalizing 
results of Ekedahl and Oort in the Siegel modular case, see~[\FOTexel].

Combining Theorem~\refn{AutDimThm} with [\Wedh], Theorem~(6.10) we arrive at the following dimension formula.

\sssection{DimForm}
{\it Corollary. --- If $\cA_{0,\delta}(w) \neq \emptyset$ then all its irreducible components have dimension equal 
to~$\ell(\dot w)$.}

\sssection{NonRedRem}
{\it Remark.\/}\enspace --- Let $\ul{Y}$ be a BT$_1$ with $(\cO/p\cO,\ast,\epsilon)$-structure associated to a 
$k$-valued point $s \in \cA_{0,\delta}(w)$. The automorphism group scheme $\bAut(\ul{Y})$ is in general (highly) 
non-reduced. In fact, [\Wedh], (2.7) shows that the tangent space of $\bAut(\ul{Y})$ at the identity is isomorphic to 
the tangent space of the deformation functor $\Def(\ul{Y})$. But by Serre-Tate and [\Wedh], (2.17) the latter is 
isomorphic to the tangent space of~$\cA_{0,\delta}$ at the point~$s$. Hence the Lie algebra of $\bAut(\ul{Y})$ has 
dimension equal to $\dim(\cA_0)$, independent of~$w$. Cf.\ Example~(\refn{AutDimExa}).
 
\ssection{Some examples}{}
\Bskip

\noindent
An immediate consequence of~\refn{DimForm} is the following.
\Cskip

\sssection{SpecStrata}
{\it Corollary. --- There is a unique Ekedahl-Oort stratum that is open in~$\cA_{0,\delta}$. There is a unique $w \in W_{\mX^0}\backslash W_{G^0}$ such that the stratum $\cA_{0,\delta}(w)$, when not empty, is $0$-dimensional.}
\Cskip

The open Ekedahl-Oort stratum, corresponding to the class of the longest element of~$W_{G^0}$, plays the role of the (generalized) ordinary locus in~$\cA_{0,\delta}$. It is studied in detail in our paper~[\STPEL]. The unique $0$-dimensional stratum (if not empty) takes the role of what in the Siegel modular case is called the superspecial stratum. (In the Siegel modular case it corresponds to the abelian varieties that are isomorphic to~$E^g$ with $E$ a supersingular elliptic curve).

Already in the simplest examples we find that for moduli points $\ul{A} = (A,\bar\lambda,\iota,\bar\eta)$ in the ``generalized'' ordinary locus, the underlying abelian variety~$A$ is not necessarily ordinary in the classical sense. A similar remark applies to the superspecial stratum. We shall illustrate this with an example.

\sssection{Exa}
We consider an example of type AL (see~\refn{BasCas}); in particular this means that there is no invariant~$\delta$ to consider.

Fix PEL data as in~\refn{PELdata} with $\cB = \tilde K$ a CM-field of degree~$2m$ over~$\mQ$. Let $K \subset \tilde K$ be the totally real subfield. Suppose $p$ is a prime number that is totally inert (unramified) in the extension $\mQ \subset K$ and that splits in the extension $K \subset \tilde K$. Set $q := p^m$. Then $B = O_\cB/pO_\cB$ is isomorphic to $\mF_q \times \mF_q$, with involution (induced by complex conjugation on~$\cB$) given by $(x,y)^\ast = (y,x)$. Fixing such an isomorphism, every BT$_1$ with $(B,\ast,-1)$-structure is of the form $\ul{Y} \times \ul{Y}^D$, where $\ul{Y}$ is a BT$_1$ with $\mF_q$-structure, and where the $(-1)$-duality $\lambda\colon Y \times Y^D \to Y^D \times Y$ is given by $(y,\eta) \mapsto (\eta,-y)$. In this way the study of BT$_1$ with $(B,\ast,-1)$-structure reduces to the study of BT$_1$ with $\mF_q$-structure (without polarization).

Via the above mechanism, a $k$-valued point of~$\cA_0$ gives rise to a BT$_1$ with $\mF_q$-structure. Its type $(d,\gf)$ is determined by the chosen PEL data. More precisely, $d = \dim_\cB(\cV)$ and the CM-type~$\gf$ is determined by the choice of the $\cG(\mR)$-conjugacy class $\cX \subset \Hom(\mS,\cG_\mR)$.

For instance, suppose $m=6$ and $d = 5$. The set $\cI$ of embeddings of $\mF_q$ into~$k = \kbar$ is a set of~$6$ elements, equipped with a natural cyclic ordering. Suppose that the function $\gf\colon \cI \to \{0,1,\ldots,5\}$ takes consecutive values $1$, $3$, $1$, $1$, $0$ and~$5$. (This is a random choice.) The BT$_1$ corresponding to the stratum of maximal dimension then corresponds to the Dieudonn\'e module given in Figure~\refn{Fig:ordExa}.
\midinsert
\centerline{\epsfbox{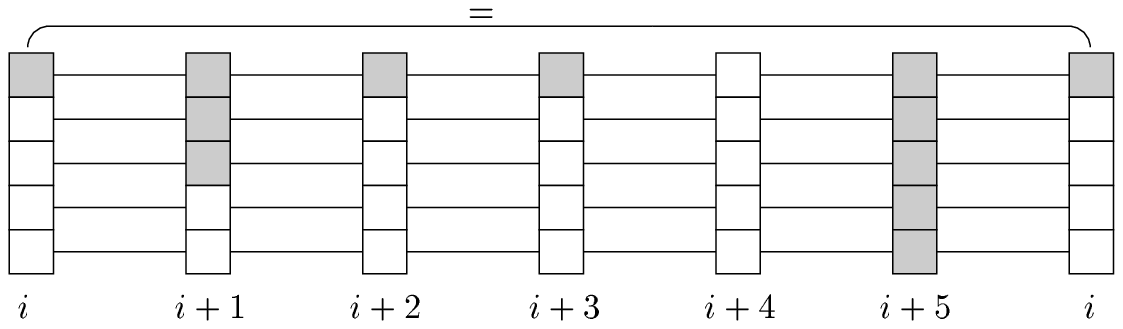}}
\figlabel{Fig:ordExa}
\endinsert

We see that in this example $\ul{Y}$ decomposes as a product, $\ul{Y} = \ul{Y}_1 \times \cdots \ul{Y}_5$, corresponding to the horizontal ``layers'' in the picture. This decomposition is not canonical but it can be shown that the coarser decomposition (referred to as the slope decomposition)
$$
\ul{Y} = (\ul{Y}_1 \times \ul{Y}_2) \times (\ul{Y}_3 \times \ul{Y}_4) \times \ul{Y}_5
$$
{\it is\/} canonical. It plays an important role in our generalization of Serre-Tate theory; we refer to~[\STPEL] for further discussion. Note, however, that $Y$, the underlying BT$_1$ without additional structure, is not ordinary in the classical sense, i.e., it is not a product of factors $\mZ/p\mZ$ and~$\mu_p$. In fact, $Y$ is a local-local group scheme.

At the other extreme, the BT$_1$ that corresponds to the $0$-dimensional stratum is given by the Dieudonn\'e module in Figure~\refn{Fig:ssExa}. The underlying group scheme~$Y$ is in this case indecomposable, quite in contrast to what happens on the superspecial stratum in the Siegel modular case.
\midinsert
\centerline{\epsfbox{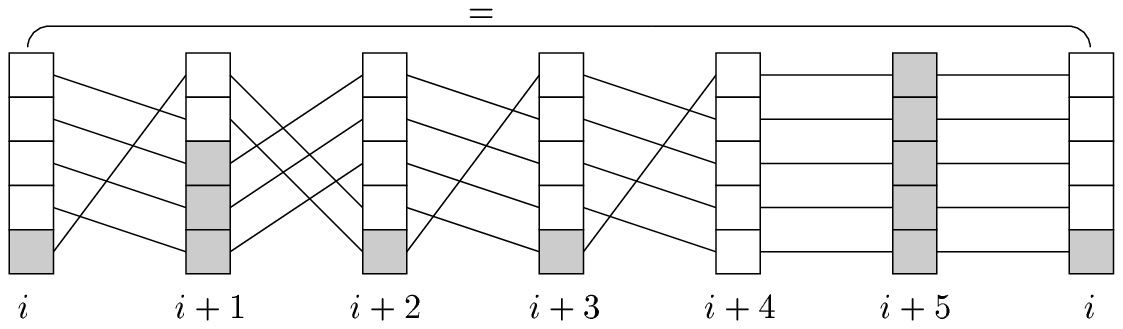}}
\figlabel{Fig:ssExa}
\endinsert

A variety of different examples is obtained by choosing different CM-types~$\gf$. Let us add that this is typical for examples of type~A (cases AU or~AL in \refn{BasCas}). In examples of type~C there is no freedom of choice for the CM-type (due to the fact that the involution is of the first kind), and we find that the BT$_1$ (without additional structure) occurring on the stratum of maximal (resp.\ minimal) dimension is of the form $\big((\mZ/p\mZ) \times \mu_p\big)^g$ (resp.\ $E[p]^g$, for $E$ a supersingular curve), as in the classical case.

\Askip\noindent
{\sectitlefont References}
\nobreak\bigskip\nobreak
{\eightpoint

\item{[\Berth]}
P.~Berthelot, {\it Th\'eorie de Dieudonn\'e sur un anneau de valuation 
parfait\/}, Ann.\ Scient.\ \'Ec.\ Norm.\ Sup.\ (4) {\bf 13} (1980), 225--268.

\item{[\Bourb]}
N.~Bourbaki, {\it Groupes et alg\`ebres de Lie, Chap.\ 4, 5, et 6\/}, Masson, 
Paris, 1981.

\item{[\deJO]}
A.J.~de Jong and F.~Oort, {\it Purity of the stratification by Newton polygons\/}, 
J.\ A.M.S.\ {\bf 13} (2000), 209--241.

\item{[\Font]}
J.-M.~Fontaine, {\it Groupes $p$-divisibles sur les corps locaux\/}, Ast\'erisque 
{\bf 47--48} (1977).  

\item{[\GorOo]}
E.~Goren and F.~Oort, {\it Stratifications of Hilbert modular varieties\/}, J.\ 
Alg.\ Geom.\ {\bf 9} (2000), 111--154.

\item{[\Illus]}
L.~Illusie, {\it D\'eformations de groupes de Barsotti-Tate (d'apr\`es 
A.~Grothendieck)\/}, in: S\'emi\-naire sur les pinceaux arithm\'etiques: la 
conjecture de Mordell (L.~Szpiro, ed.), Ast\'erisque {\bf 127} (1985), 151--198.

\item{[\KottIsoc]}
R.E.~Kottwitz, {\it Isocrystals with additional structure\/}, Compos.\ Math.\ 
{\bf 56} (1985), 201--220; {\it ---~\Romno 2\/}, Compos.\ Math.\ {\bf 109} (1997), 
255--339.

\item{[\Kott]}
R.E.~Kottwitz, {\it Points on some Shimura varieties over finite fields\/}, 
J.A.M.S.\ {\bf 5} (1992), 373--444.

\item{[\Kraft]}
H.~Kraft, {\it Kommutative algebraische $p$-Gruppen (mit Anwendungen auf 
$p$-divisible Gruppen und abelsche Variet\"aten)\/}, manuscript, Univ.\ Bonn, 
Sept.\ 1975, 86 pp. (Unpublished)

\item{[\GSAS]}
B.J.J.~Moonen, {\it Group schemes with additional structures and Weyl group 
cosets\/}, in: Moduli of Abelian Varieties (C.~Faber, G.~van der Geer, F.~Oort, 
eds.), Progr.\ Math.\ {\bf 195}, Birkh\"auser, Basel, 2001, pp.\ 255--298.

\item{[\STPEL]}
B.J.J.~Moonen, {\it Serre-Tate theory for moduli spaces of PEL type\/}, 
math.AG/0203288.

\item{[\FOAnn]}
F.~Oort, {\it Newton polygons and formal groups: conjectures by Manin and 
Grothendieck\/}, Ann.\ Math.\ {\bf 152} (2000), 183--206.

\item{[\FOTexel]}
F.~Oort, {\it A stratificiation of a moduli space of abelian varieties\/}, in:
Moduli of Abelian Varieties (C.~Faber, G.~van der Geer, F.~Oort, eds.), Progr.\ 
Math.\ {\bf 195}, Birkh\"auser, Basel, 2001, pp.\ 345--416.

\item{[\FOTexelB]}
F.~Oort, {\it Newton polygon strata in the moduli space of abelian varieties\/}, in:
Moduli of Abelian Varieties (C.~Faber, G.~van der Geer, F.~Oort, eds.), Progr.\ 
Math.\ {\bf 195}, Birkh\"auser, Basel, 2001, pp.\ 417--440.

\item{[\RapBourb]}
M.~Rapoport, {\it On the Newton stratification\/}, S\'em.\ Bourbaki, Mars 2002.

\item{[\RR]}
M.~Rapoport and M.~Richartz, {\it On the classification and specialization of 
$F$-isocrystals with additional structure\/}, Compos.\ Math.\ {\bf 103} (1996), 
153--181.

\item{[\Wedh]}
T.~Wedhorn, {\it The dimension of Oort strata of Shimura varieties of PEL-type\/}, 
in: Moduli of Abelian Varieties (C.~Faber, G.~van der Geer, F.~Oort, eds.), Progr.\ 
Math.\ {\bf 195}, Birkh\"auser, Basel, 2001, pp.\ 441--471.
\Bskip 

\noindent
Ben Moonen, University of Amsterdam, Korteweg-de Vries Institute for Mathematics, 
Plantage Muidergracht~24, 1018~TV Amsterdam, The Netherlands. 
Email: {\tt bmoonen@science.uva.nl}\par}
\bye